\documentclass[reqno]{amsart}

\usepackage[latin1]{inputenc}
\usepackage{hyperref}
\usepackage{amssymb,dsfont}

\newcommand{\restr}[1]{_{\rule[-.41ex]{.03em}{2ex}#1}}
\newcommand{\intn}[2]{\ensuremath{\{  #1  ,\ldots,  #2 \}}}
\newcommand{\gene}[1]{\ensuremath{\langle #1 \rangle}}

\def\grint{\mathrm{G}}

\def\Id{\mathrm{Id}}

\def\R{\mathbb{R}}

\def\N{\mathbb{N}}
\def\Z{\mathbb{Z}}
\def\T{\mathbb{T}}

\DeclareMathAlphabet{\mathpzc}{OT1}{pzc}{m}{it}
\def\dist{\mathpzc{d}}

\def\hau{\mathcal{H}}
\def\leb{\mathcal{L}}
\def\netm{\mathcal{M}}

\def\prob{\mathbb{P}}
\def\esp{\mathbb{E}}
\def\un{\mathds{1}}
\def\dd{\mathrm{d}}

\def\I{\mathcal{I}}

\def\netm{\mathcal{M}}
\newcommand{\diam}[1]{\ensuremath{| #1 |}}

\def\eps{\varepsilon}
\def\ph{\varphi}

\def\jauge{\mathfrak{D}}
\def\jaugew{\mathcal{D}}
\def\modc{\mathfrak{W}}
\def\HH{\mathbb{H}}

\def\Mu{\mathrm{M}}
\def\Ncal{\mathcal{N}}

\theoremstyle{plain}
\newtheorem{thm}{Theorem}
\newtheorem{prp}{Proposition}
\newtheorem{lem}[prp]{Lemma}

\theoremstyle{definition}
\newtheorem*{df}{Definition}

\theoremstyle{remark}
\newtheorem*{rems}{Remarks}
\newtheorem*{rem}{Remark}

\title{Singularity sets of Lévy processes}
\author{Arnaud Durand}
\subjclass[2000]{Primary 60D05; Secondary 60G51, 60G17, 26A15, 28A80}
\address{Laboratoire d'Analyse et de Mathématiques Appliquées, Université Paris XII, 61 av. du Général de Gaulle, 94010 Créteil Cedex, France.}
\email{a.durand@univ-paris12.fr}

\begin{document}

\begin{abstract}
We completely describe the size and large intersection properties of the H\"older singularity sets of Lévy processes. We also study the set of times at which a given function cannot be a modulus of continuity of a Lévy process. The H\"older singularity sets of the sample paths of certain random wavelet series are investigated as well.
\end{abstract}

\maketitle

\section{Introduction}\label{grintlevyintro}

Let us consider a $d$-dimensional Lévy process $X=(X_{t})_{t\in [0,\infty)}$ whose Lévy measure has infinite total mass. A remarkable property enjoyed by almost all sample paths of $X$ is that they are {\em multifractal functions}. This property, which was proven by S.~Jaffard~\cite{Jaffard:1999fg}, implies that the regularity of $X$ fluctuates so erratically as time passes that with probability one, the random sets
\begin{equation}\label{defEh}
E_{h}=\left\{ t\in [0,\infty) \:\bigl|\: h_{X}(t)=h \right\}
\end{equation}
are nonempty for all $h$ in some subinterval of $[0,\infty]$ which is not reduced to a point. Here, $h_{X}(t)$ denotes the {\em H\"older exponent} of the process $X$ at time $t$. It measures the regularity of $X$ at $t$ and is defined as the supremum of all $h>0$ such that there are a real $c>0$ and a $d$-tuple $P$ of polynomials enjoying
\[
\|X_{t'}-P(t'-t)\|\leq c\, |t'-t|^h
\]
for any nonnegative real $t'$ in a neighborhood of $t$, see~\cite{Jaffard:2004fh}. S.~Jaffard actually established a more precise result since he determined the {\em spectrum of singularities} $h\mapsto\dim E_{h}$ of almost every sample path of $X$ (see Theorem~\ref{speclevythm} below). Here, $\dim$ stands for Hausdorff dimension. He thus gave a first description of the size properties of the sets $E_{h}$.

In this paper, we provide a finer description of these size properties by computing the Hausdorff $g$-measure of the sets $E_{h}\cap V$ for any gauge function $g$ and every open set $V$. On top of that, we show that certain random sets related to the sets $E_{h}$ enjoy a remarkable geometrical property which was introduced by K.~Falconer~\cite{Falconer:1994hx}. To be specific, we establish that the sets
\begin{equation}\label{deftildeEh}
\tilde E_{h}=\left\{ t\in [0,\infty)\backslash S \:\bigl|\: h_{X}(t)\leq h \right\},
\end{equation}
where $S$ denotes the set of jump times of the Lévy process $X$, are almost surely {\em sets with large intersection}. This means in particular that they are locally everywhere of the same size, in the sense that the Hausdorff dimension of $\tilde E_{h}\cap V$ does not depend on the choice of the nonempty open subset $V$ of $(0,\infty)$. This also implies that the size properties of the sets $\tilde E_{h}$ are not altered by taking countable intersections. Indeed, the Hausdorff dimension of the intersection of countably many sets with large intersection is equal to the infimum of their Hausdorff dimensions. This property is rather counterintuitive in view of the fact that the intersection of two subsets of $\R$ of Hausdorff dimensions $s_{1}$ and $s_{2}$ respectively is usually expected to be $s_{1}+s_{2}-1$, see~\cite[Chapter~8]{Falconer:2003oj} for precise statements. The occurrence of sets with large intersection in the theory of Diophantine approximation and that of dynamical systems was pointed out by many authors, see~\cite{Durand:2006uq,Durand:2007uq,Falconer:1994hx} and the references therein. Their use in multifractal analysis of stochastic processes is more novel and was introduced by J.-M.~Aubry and S.~Jaffard~\cite{Aubry:2002up} in order to determine the law of the spectrum of singularities of a model of random wavelet series. Our results indicate that sets with large intersection also arise in the study of more common processes like Lévy processes.

We also study the set of times at which a given function cannot be a modulus of continuity of a Lévy process. To be specific, we show that this set almost surely contains a set with large intersection and we give a sufficient condition on a gauge function $g$ to ensure that this set almost surely has maximal Hausdorff $g$-measure in every open subset of $(0,\infty)$.

The methods that we develop hereunder to investigate the size and large intersection properties of the singularity sets $E_{h}$ and $\tilde E_{h}$ of Lévy processes can also be applied to a model of random wavelet series which generalizes that previously introduced by S.~Jaffard in~\cite{Jaffard:2000mk}. We thus obtain analogous results in that context, see Section~\ref{grintlevylac}.

\section{Statement of the results}

We begin by recalling some basic properties of Lévy processes. Remember that $X=(X_{t})_{t\in [0,\infty)}$ denotes a Lévy process valued in $\R^d$. This means that $X$ has stationary independent increments, its sample paths are right-continuous with left limits and $X_{0}$ vanishes. The characteristic exponent $\psi$ is defined by $\esp[e^{i \lambda\cdot X_{t}}]=e^{-t\psi(\lambda)}$ for all $t\in [0,\infty)$ and $\lambda\in\R^d$ and is given by the Lévy-Khintchine formula, that is,
\[
\psi(\lambda)=i a\cdot\lambda+\frac{1}{2}q(\lambda)+\int_{\R^d} \left( 1-e^{i \lambda\cdot x}+i \lambda\cdot x\ \un_{\{\|x\|<1\}} \right) \pi(\dd x)
\]
where $a\in\R^d$, $q$ is a positive semidefinite quadratic form on $\R^d$ and $\pi$ is a measure on $\R^d\backslash\{0\}$ which enjoys $\int (1\wedge\|x\|^2)\pi(\dd x)<\infty$ and is called the Lévy measure of $X$, see e.g.~\cite[Theorem~8.1]{Sato:1999fk}. Note that the local regularity properties of $X$ are trivial if the total mass of $\pi$ is finite, because in this case $X$ is the superposition of a compound Poisson process with drift and a Brownian motion, whose H\"older exponent is~$1/2$ everywhere with probability one, see e.g.~\cite[Chapter~16]{Kahane:1985gc} or~\cite[Chapter~2]{Karatzas:1991tv}. Therefore, our results are nontrivial only if the Lévy measure $\pi$ has infinite total mass. Moreover, let $\Delta X_{t}=X_{t}-X_{t^{-}}$ for each time $t\in (0,\infty)$ and let $S$ denote the set of all times $t$ such that $\Delta X_{t}\neq 0$. Thus, $S$ is the set of all jump times of the Lévy process $X$.

\subsection{Size properties of the singularity sets}

We shall give an exhaustive description of the size properties of the sets $E_{h}$ and $\tilde E_{h}$ defined by~(\ref{defEh}) and~(\ref{deftildeEh}) respectively for every $h\in [0,\infty]$.

A typical way to describe the size properties of a subset of $\R^d$ is to compute its Hausdorff $g$-measure in every open subset of $\R^d$, for every gauge function $g$. Recall that a gauge function is a nondecreasing function $g$ defined on $[0,\eps]$ for some $\eps>0$ and such that $\lim_{0^+}g=g(0)=0$. The Hausdorff $g$-measure of a set $F\subseteq\R^d$ is defined by
\[
\hau^g(F)=\lim_{\delta\downarrow 0}\uparrow \inf_{F\subseteq\bigcup_{p} U_{p} \atop \diam{U_{p}}<\delta} \sum_{p=1}^\infty g(\diam{U_{p}}).
\]
The infimum is taken over all sequences $(U_{p})_{p\in\N}$ of subsets of $\R^d$ with $F\subseteq\bigcup_{p} U_{p}$ and $\diam{U_{p}}<\delta$ for all $p\in\N$, where $\diam{\cdot}$ denotes diameter. As stated in~\cite{Rogers:1970wb}, $\hau^g$ is a Borel measure on $\R^d$. Moreover, if $g(r)/r^d$ tends to infinity as $r\to 0$, then any subset of $\R^d$ with nonvanishing Lebesgue measure has infinite Hausdorff $g$-measure. Otherwise, $\hau^g$ is a translation invariant Borel measure which is finite on compacts, so that it coincides up to a multiplicative constant with the Lebesgue measure on the Borel subsets of $\R^d$.

Restricting to the gauge functions $\Id^s$, where $\Id$ denotes the identity function, leads to the notion of Hausdorff dimension. Specifically, the Hausdorff dimension of a nonempty set $F\subseteq\R^d$ is defined by
\[
\dim F=\sup\{ s\in (0,d) \:|\: \hau^{\Id^s}(F)=\infty \}=\inf\{ s\in (0,d) \:|\: \hau^{\Id^s}(F)=0 \},
\]
thus giving a partial description of its size properties, see~\cite{Falconer:2003oj}. In the above formula, we adopt the convention that $\sup\emptyset=0$ and $\inf\emptyset=d$.

S.~Jaffard~\cite{Jaffard:1999fg} computed the Hausdorff dimension of the sets $E_{h}$. In order to state his result, let us introduce some further notations. Let $\sigma$ be the image measure of $\pi$ by $x\mapsto \|x\|$ and let $c_{j}=\sigma((2^{-j-1},2^{-j}])$ for every integer $j\geq 0$. We shall often assume that the following condition holds:
\begin{equation}\label{condcjlevy}
\sum_{j=0}^\infty 2^{-j}\sqrt{c_{j}\log(1+c_{j})}<\infty.
\end{equation}
Note that~(\ref{condcjlevy}) is satisfied as soon as the exponent
\[
\beta=\inf\left\{ \gamma\in [0,\infty) \:\biggl|\: \int_{0}^1 r^\gamma \sigma(\dd r)<\infty \right\}
\]
is less than $2$. This exponent, which lies in $[0,2]$, was introduced by R.~Blumenthal and R.~Getoor~\cite{Blumenthal:1961kx} and has an effect on the pointwise regularity of the Lévy process $X$. First, W.~Pruitt~\cite{Pruitt:1981vn} proved that $h_{X}(0)=1/\beta$ almost surely if $q=0$. Thus the H\"older exponent of $X$ is $1/\beta$ almost everywhere with probability one. Then, let
\[
\beta'=\begin{cases}
\beta & \text{if }q=0 \\
2 & \text{if }q\neq 0
\end{cases}
\]
with the convention that $1/\beta'=\infty$ if $\beta'=0$. As shown by the following result of S.~Jaffard, the exponents $\beta$ and $\beta'$ govern the law of the Hausdorff dimension of the sets $E_{h}$.

\begin{thm}[S.~Jaffard]\label{speclevythm}
Suppose that~(\ref{condcjlevy}) holds.\begin{enumerate}
\item With probability one, for every $h\in (1/\beta',\infty]$, $E_{h}=\emptyset$.
\item With probability one, $E_{1/\beta'}$ has full Lebesgue measure in $[0,\infty)$.
\item If $\beta>0$, with probability one, for every $h\in [0,1/\beta')$, $\dim E_{h}=\beta h$.
\item If $\beta=0$ and $\pi(\R^d)=\infty$, for any fixed $h\in [0,1/\beta')$, $\dim E_{h}=0$ with probability one.
\end{enumerate}\end{thm}

We shall refine this result by computing $\hau^g(E_{h}\cap V)$ for every gauge function $g\in\jauge$ and every open subset $V$ of $\R$ and we shall also determine $\hau^g(\tilde E_{h}\cap V)$, thereby providing a full description of the size properties of the sets $E_{h}$ and $\tilde E_{h}$. We can actually restrict our attention to the case where $V$ is included in $(0,\infty)$. Indeed, for any open subset $V$ of $\R$, we have $\hau^g(E_{h}\cap V)=\hau^g(E_{h}\cap V_{+})$ and $\hau^g(\tilde E_{h}\cap V)=\hau^g(\tilde E_{h}\cap V_{+})$, where $V_{+}$ denotes the open subset $V\cap (0,\infty)$ of $(0,\infty)$. In addition, the size properties of $E_{h}$ and $\tilde E_{h}$ are trivial when $h\in [1/\beta',\infty]$. More precisely, Theorem~\ref{speclevythm} ensures that with probability one, $E_{h}=\emptyset$ and $\tilde E_{h}=\tilde E_{1/\beta'}$ for all $h\in (1/\beta',\infty]$. Furthermore, along with the fact that $S$ is almost surely countable, Theorem~\ref{speclevythm} implies that the sets $E_{1/\beta'}$ and $\tilde E_{1/\beta'}$ almost surely have full Lebesgue measure in $[0,\infty)$, so that with probability one, for every $g\in\jauge$ and every open subset $V$ of $(0,\infty)$,
\[
\hau^g(E_{1/\beta'}\cap V)=\hau^g(\tilde E_{1/\beta'}\cap V)=\hau^g(V).
\]

For those reasons, the only case of interest is that in which $V\subseteq (0,\infty)$ and $h\in [0,1/\beta')$. In this case, the values of $\hau^g(E_{h}\cap V)$ and $\hau^g(\tilde E_{h}\cap V)$ for $g\in\jauge$ are given by Theorems~\ref{grintisoholdlevy1} and~\ref{grintisoholdlevy2} below. In order to state these results, we need to introduce some additional notations. Let $\jauge_{d}$ denote the set of all gauge functions $g\in\jauge$ such that $r\mapsto g(r)/r^d$ is positive and nonincreasing on $(0,\eps]$ for some $\eps>0$. One easily checks that any function in $\jauge_{d}$ is continuous in a neighborhood of zero. For any $g\in\jauge$, let
\begin{equation}\label{defjaugeg1}
g_{d}:r\mapsto r^d\inf_{\rho\in (0,r]}\frac{g(\rho)}{\rho^d}.
\end{equation}
Proposition~2 in~\cite{Durand:2007uq} then ensures that $g_{d}\in\jauge_{d}\cup\{0\}$ for any $g\in\jauge$ and that there exists a real number $\kappa\geq 1$ such that
\begin{equation}\label{prpjauge1}
\forall g\in\jauge \quad \forall F\subset\R^d \qquad \hau^{g_{d}}(F)\leq\hau^g(F)\leq\kappa\,\hau^{g_{d}}(F).
\end{equation}
Here, $\kappa$ only depends on the dimension $d$ and the norm $\R^d$ is endowed with. This means that the Hausdorff measures built using the gauge functions $g$ and $g_{d}$ are comparable. Therefore, to study size properties of subsets of $\R^d$, one often begins by considering gauge functions in $\jauge_{d}$.

Moreover, for any gauge function $g$ in $\jauge$, let
\[
h_{g}=\inf\left\{ h\in (0,\infty) \:\biggl|\: \int_{0^{+}} g_{1}(r^{1/h}) \sigma(\dd r)=\infty \right\}
\]
where $\int_{0^+}$ denotes the integral on any interval $(0,\eps]$ with $\eps>0$ on which the integrand is continuous. In the above formula, we adopt the usual convention that $\inf\emptyset=\infty$. It is straightforward to prove that $h_{g}=\infty$ if $g_{1}=0$ and that $h_{g}\leq 1/\beta$ otherwise. Besides, $h_{\Id^s}=s/\beta$ for any $s\in (0,1]$.

\begin{thm}\label{grintisoholdlevy1}
Assume that~(\ref{condcjlevy}) holds and let $g\in\jauge$. With probability one, for every real $h\in [0,1/\beta')$ and every nonempty open subset $V$ of $(0,\infty)$,
\[
\hau^g(E_{h}\cap V)=\begin{cases}
0 & \text{if }h<h_{g} \\
\infty & \text{if }h=h_{g}
\end{cases}
\qquad\text{and}\qquad
\hau^g(\tilde E_{h}\cap V)=\begin{cases}
0 & \text{if }h<h_{g} \\
\infty & \text{if }h\geq h_{g}
\end{cases}.
\]
\end{thm}

\begin{rems}
This result shows that for any $g\in\jauge$ with $h_{g}=0$, the set $E_{0}$ almost surely has infinite Hausdorff $g$-measure in every nonempty open subset $V$ of $(0,\infty)$. Such a function $g$ exists if and only if the Lévy measure $\pi$ has infinite total mass. Thus, if its Lévy measure has infinite total mass, there are uncountably many times at which a Lévy process has a vanishing H\"older exponent and does not jump.

Let $h\in (0,1/\beta')$. Theorem~\ref{grintisoholdlevy1} ensures that $E_{h}$ almost surely has infinite Hausdorff $g$-measure in every nonempty open subset of $(0,\infty)$ for any fixed gauge function $g\in\jauge$ with $h_{g}=h$. If $\beta>0$, such a function always exists: take $g=\Id^{\beta h}$ for example. This need not be true if $\beta=0$. Indeed, some Lévy measures $\pi$ yield $h_{g}\in\{0,\infty\}$ for all $g\in\jauge$. Consider e.g.~$\pi=\sum_{j=1}^\infty\delta_{1/j}/j$ in dimension $d=1$.
\end{rems}

Observe that the event of probability one on which the statement of Theorem~\ref{grintisoholdlevy1} holds depends on the gauge function $g$. One can gain uniformity in $g$ by assuming that $\beta>0$ and restricting oneself to the following collection of gauge functions:
\[
\jaugew=\left\{ g\in\jauge \:\biggl|\: \liminf_{r\to 0}\frac{\log g_{1}(r)}{\log r}=\limsup_{r\to 0}\frac{\log g_{1}(r)}{\log r} \right\}\subseteq\jauge.
\]
Note that for any $s\in (0,1]$ the function $\Id^s$ belongs to $\jaugew$.

\begin{thm}\label{grintisoholdlevy2}
Assume that $\beta>0$ and that~(\ref{condcjlevy}) holds. With probability one, for every $g\in\jaugew$, every $h\in [0,1/\beta')$ and every nonempty open subset $V$ of $(0,\infty)$,
\[
\hau^g(E_{h}\cap V)=\hau^g(\tilde E_{h}\cap V)=\begin{cases}
0 & \text{if }h<h_{g} \\
\infty & \text{if }h\geq h_{g}
\end{cases}.
\]
\end{thm}

If $\beta>0$ and~(\ref{condcjlevy}) holds, it obviously follows from Theorem~\ref{grintisoholdlevy2} that with probability one, for every $h\in [0,1/\beta')$ and $s\in (0,1]$,
\[
\hau^{\Id^s}(E_{h})=\begin{cases}
0 & \text{if }s>\beta h \\
\infty & \text{if }s\leq \beta h
\end{cases}
\]
so that $\dim E_{h}=\beta h$. Thus Theorem~\ref{grintisoholdlevy2} yields the part of Theorem~\ref{speclevythm} that concerns the sets $E_{h}$ for $h\in [0,1/\beta')$.

\subsection{Large intersection properties of the singularity sets}

We shall show that the sets $\tilde E_{h}$ are sets with large intersection, in the sense that they belong to certain classes $\grint^g(V)$ of subsets of $\R$. These classes were introduced in~\cite{Durand:2007uq} in order to generalize the original classes of sets with large intersection of K.~Falconer~\cite{Falconer:1994hx}.

The classes $\grint^g(V)$ are defined as follows in general dimension $d$. For any $g\in\jauge_{d}$, let $\eps_{g}$ denote the supremum of all $\eps\in (0,1]$ such that $g$ is nondecreasing on $[0,\eps]$ and $r\mapsto g(r)/r^d$ is nonincreasing on $(0,\eps]$. Moreover, given an integer $c\geqslant 2$, let $\Lambda_{c,g}$ be the set of all $c$-adic cubes of diameter less than $\eps_{g}$, that is, sets of the form $\lambda=c^{-j}(k+[0,1)^d)$ with $j\in\Z$, $k\in\Z^d$ and $\diam{\lambda}<\eps_{g}$. The outer net measure associated with $g\in\jauge_{d}$ is defined by
\begin{equation}\label{defnetmh}
\forall F\subseteq\R^d \qquad \netm^g_{\infty}(F) = \inf_{(\lambda_{p})_{p\in\N}} \sum_{p=1}^\infty g(\diam{\lambda_{p}}).
\end{equation}
The infimum is taken over all sequences $(\lambda_{p})_{p\in\N}$ in $\Lambda_{c,g}\cup\{\emptyset\}$ enjoying $F\subseteq\bigcup_{p}\lambda_{p}$. The outer measure $\netm^g_{\infty}$ is related to $\hau^g$, see~\cite[Theorem~49]{Rogers:1970wb}. In particular, if a subset $F$ of $\R^d$ enjoys $\netm^g_{\infty}(F)>0$ then $\hau^g(F)>0$. In addition, for $\overline{g},g\in\jauge_{d}$, let us write $\overline{g}\prec g$ if $\overline{g}/g$ monotonically tends to infinity at zero. We can now define the classes $\grint^g(V)$. Recall that a $G_{\delta}$-set is one that may be expressed as a countable intersection of open sets.

\begin{df}
Let $g\in\jauge_{d}$ and let $V$ be a nonempty open subset of $\R^d$. The class $\grint^g(V)$ of subsets of $\R^d$ with large intersection in $V$ with respect to $g$ is the collection of all $G_{\delta}$-subsets $F$ of $\R^d$ such that $\netm^{\overline{g}}_{\infty}(F\cap U)=\netm^{\overline{g}}_{\infty}(U)$ for every $\overline{g}\in\jauge_{d}$ enjoying $\overline{g}\prec g$ and every open set $U\subseteq V$.
\end{df}

\begin{rem}
It is proven in~\cite{Durand:2007uq} that the class $\grint^g(V)$ depends on the choice of neither the integer $c$ nor the norm $\R^d$ is endowed with, even if they affect the construction of $\netm^{\overline{g}}_{\infty}$ for any $\overline{g}\in\jauge_{d}$ with $\overline{g}\prec g$.
\end{rem}

The class $\grint^g(V)$ enjoys several remarkable properties, among which the following are the most important, see~\cite{Durand:2007uq}.

\begin{thm}\label{GRINTSTABLE}
Let $g\in\jauge_{d}$ and let $V$ be a nonempty open subset of $\R^d$. Then\begin{enumerate}
\item the class $\grint^g(V)$ is closed under countable intersections;
\item the set $f^{-1}(F)$ belongs to $\grint^g(V)$ for every bi-Lipschitz mapping $f:V\to\R^d$ and every set $F\in\grint^g(f(V))$;
\item every set $F\in\grint^g(V)$ enjoys $\hau^{\overline{g}}(F)=\infty$ for every $\overline{g}\in\jauge_{d}$ with $\overline{g}\prec g$ and in particular $\dim F\geq s_{g}=\sup\{ s\in (0,d) \:|\: \Id^s\prec g \}$;
\item every $G_{\delta}$-subset of $\R^d$ of full Lebesgue measure in $V$ belongs to $\grint^g(V)$.
\end{enumerate}
\end{thm}

For $h\in [0,\infty]$, we shall determine for which gauge functions $g\in\jauge_{1}$ and which nonempty open subsets $V$ of $\R$ the set $\tilde E_{h}$ belongs to the class $\grint^g(V)$. It is easy to check that $\tilde E_{h}$, being included in $[0,\infty)$, cannot belong to $\grint^g(V)$ if $V\not\subseteq (0,\infty)$. Furthermore, Theorem~\ref{speclevythm} and the observation that $S$ is almost surely countable imply that $\tilde E_{1/\beta'}$ has full Lebesgue measure in $[0,\infty)$ with probability one. It follows from Theorem~\ref{GRINTSTABLE} that with probability one, for every $g\in\jauge_{1}$, every $h\in [1/\beta',\infty]$ and every nonempty open subset $V$ of $(0,\infty)$, the set $\tilde E_{h}$ belongs to $\grint^g(V)$.

Therefore, the only case of interest is that in which $V\subseteq (0,\infty)$ and $h\in [0,1/\beta')$. In this case, the gauge functions $g$ and the open sets $V$ such that $\tilde E_{h}\in\grint^g(V)$ are given by the next theorem. In its statement, $\jaugew_{1}$ denotes the collection $\jaugew\cap\jauge_{1}$ of gauge functions.

\begin{thm}\label{grintisoholdlevy3}
Assume that~(\ref{condcjlevy}) holds.
\begin{enumerate}
\item\label{grintisoholdlevy31} Let $g\in\jauge_{1}$. With probability one, for every $h\in [0,1/\beta')$ and every nonempty open subset $V$ of $(0,\infty)$,
\[
\tilde E_{h}\in\grint^g(V) \qquad\Longleftrightarrow\qquad h\geq h_{g}.
\]
\item\label{grintisoholdlevy32} Suppose that $\beta>0$. With probability one, for every $g\in\jaugew_{1}$, every $h\in [0,1/\beta')$ and every nonempty open subset $V$ of $(0,\infty)$,
\[
\tilde E_{h}\in\grint^g(V) \qquad\Longleftrightarrow\qquad h\geq h_{g}.
\]
\end{enumerate}
\end{thm}

Let us mention that Theorems~\ref{grintisoholdlevy1},~\ref{grintisoholdlevy2} and~\ref{grintisoholdlevy3} can be seen as corollaries of a more general result, namely, Theorem~\ref{grintisoholdlevy}, which is stated and proven in Section~\ref{grintlevyhold}.

\subsection{Pointwise moduli of continuity}

We also study the size and large intersection properties of the set of all times $t\in [0,\infty)$ at which a given function cannot be a modulus of continuity of the Lévy process $X$. Let $\modc$ be the set of all continuous increasing functions $w$ defined on $[0,\eps]$ for some $\eps>0$ and such that $w(0)=0$ and
\[
1<\liminf_{\delta\to 0} \frac{w(2\delta)}{w(\delta)} \leq \limsup_{\delta\to 0} \frac{w(2\delta)}{w(\delta)} <\infty.
\]
For any $h\in (0,\infty)$, one readily verifies that the function $\delta\mapsto\delta^h$ belongs to $\modc$. Let $t\in [0,\infty)$. A function $w\in\modc$ is a modulus of continuity of $X$ at $t$ if there are a real $c>0$ and a $d$-tuple $P$ of polynomials enjoying
\[
\| X_{t'}-P(t' - t) \| \leq c\, w(| t' - t |)
\]
for any nonnegative real $t'$ in a neighborhood of $t$. Given a function $w\in\modc$, let $F_{w}$ denote the set of all times $t\in [0,\infty)\backslash S$ at which $w$ is not a modulus of continuity of $X$. The following result is established in Section~\ref{grintlevymodc}.

\begin{thm}\label{grintmodclevy}
Let $g\in\jauge$ and $w\in\modc$ with $\int_{0^{+}} g_{1}(w^{-1}(r))\sigma(\dd r)=\infty$. Then $F_{w}$ almost surely contains a set of the class $\grint^{g_{1}}((0,\infty))$. Furthermore, with probability one, for every open subset $V$ of $(0,\infty)$, we have $\hau^g(F_{w}\cap V)=\hau^g(V)$.
\end{thm}

\begin{rem}
The part of the statement of Theorem~\ref{grintmodclevy} that concerns the size properties of $F_{w}$ is a convenient way to recast two results. Let $g\in\jauge$ and $w\in\modc$ with $\int_{0^{+}} g_{1}(w^{-1}(r))\sigma(\dd r)=\infty$. On the one hand, if $g(r)/r$ tends to infinity as $r\to 0$, the set $F_{w}$ almost surely has infinite Hausdorff $g$-measure in every nonempty open subset of $(0,\infty)$. On the other hand, if $g(r)/r$ does not tend to infinity at zero, $F_{w}$ actually contains a Borel set of full Lebesgue measure in $[0,\infty)$ with probability one. The fact that $F_{w}$ almost surely has maximal Hausdorff $g$-measure in every open subset of $(0,\infty)$ then follows from the observation that $\hau^g$ coincides up to a multiplicative constant with the Lebesgue measure on the Borel subsets of $\R$. We refer to Section~\ref{grintlevymodc} for details.
\end{rem}

\section{Singularity sets}\label{grintlevyhold}

In this section, we establish the following result and we explain how it leads to Theorems~\ref{grintisoholdlevy1},~\ref{grintisoholdlevy2} and~\ref{grintisoholdlevy3}. Note that $u\circ v\in\jauge_{1}$ for every $u\in\jauge_{1}$ and every $v\in\jaugew_{1}$.

\begin{thm}\label{grintisoholdlevy}
Assume that~(\ref{condcjlevy}) holds and let $u\in\jauge_{1}$ with $h_{u}<\infty$. With probability one, for every $v\in\jaugew_{1}$, every $h\in [0,1/\beta')$ and every nonempty open $V\subseteq (0,\infty)$,
\[
\hau^{u\circ v}(E_{h}\cap V)=\begin{cases}
0 & \text{if }h<h_{u\circ v} \\
\infty & \text{if }h_{u\circ v}\leq h\leq h_{u}
\end{cases}
\,\text{and}\quad
\hau^{u\circ v}(\tilde E_{h}\cap V)=\begin{cases}
0 & \text{if }h<h_{u\circ v} \\
\infty & \text{if }h\geq h_{u\circ v}
\end{cases}
\]
and the set $\tilde E_{h}$ belongs to the class $\grint^{u\circ v}(V)$ if and only if $h\geq h_{u\circ v}$.
\end{thm}

The section is organized as follows. We begin by recalling some results obtained by S.~Jaffard in~\cite{Jaffard:1999fg} and by proving several preliminary lemmas. We then establish Theorem~\ref{grintisoholdlevy}. Theorems~\ref{grintisoholdlevy1},~\ref{grintisoholdlevy2} and~\ref{grintisoholdlevy3} are proven at the end of this section.

\subsection{Preliminary results}

Let us recall a basic property concerning the jumps of the Lévy process $X$. Let $\HH=(0,\infty)\times (\R^d\backslash\{0\})$ and, for any Borel subset $B$ of $\HH$, let $J(B)$ be the number of times $t\in (0,\infty)$ enjoying $(t,\Delta X_{t})\in B$. Then $J$ is a Poisson random measure with intensity $\leb^1\restr{(0,\infty)}\otimes\pi$, where $\leb^1$ denotes the Lebesgue measure on $\R$, see~\cite[Theorem~19.2]{Sato:1999fk}. It follows in particular that the set $S$ of jump times of $X$ is almost surely countable.

The H\"older exponent of the Lévy process $X$ at a given time $t\in [0,\infty)$ depends on the accuracy with which its jump times approach $t$. More precisely, let $S_{1}$ denote the set of all $s\in S$ such that $\|\Delta X_{s}\|\leq 1$ and, for every continuous nondecreasing function $\ph:[0,\infty)\to\R$ enjoying $\ph(0)=0$, let
\[
L_{\ph}=\left\{ t\in [0,\infty) \:\bigl|\: |t-s|<\ph(\|\Delta X_{s}\|)\text{ for infinitely many }s\in S_{1} \right\}.
\]
It is easy to check that $\alpha\mapsto L_{\Id^{1/\alpha}}$ is nondecreasing. The following proposition recasts both Lemma~2 and Proposition~1 in~\cite{Jaffard:1999fg}.

\begin{prp}[S.~Jaffard]\label{locholdlevy}
Assume that~(\ref{condcjlevy}) holds. With probability one, for every real number $h\in [0,1/\beta')$,
\[
\tilde E_{h}=\Biggl(\bigcap_{h<\alpha\leq 1/\beta} L_{\Id^{1/\alpha}}\Biggr)\backslash S \qquad\text{and}\qquad
E_{h}\backslash S=\tilde E_{h}\backslash\bigcup_{0<\alpha<h} L_{\Id^{1/\alpha}}.
\]
\end{prp}

Owing to Proposition~\ref{locholdlevy} and the fact that $S$ is almost surely countable, in order to establish Theorem~\ref{grintisoholdlevy}, it suffices to determine the size and large intersection properties of $L_{\Id^{1/\alpha}}$ for any $\alpha\in (0,\infty)$. The next lemma is a first step towards this goal.

\begin{lem}\label{sizelph}
Consider a continuous nondecreasing function $\ph:[0,\infty)\to\R$ enjoying $\ph(0)=0$. Then
\begin{equation*}\begin{split}
& \int_{0^{+}} \ph(r)\sigma(\dd r)<\infty \quad\Longrightarrow\quad \text{a.s.} \quad \forall n\in\N \quad \lim_{\eps\downarrow 0}\downarrow \int_{0\leq t\leq n \atop \|x\|\leq\eps} \ph(\|x\|) J(\dd t,\dd x)=0 \\
\text{and}\ & \int_{0^{+}} \ph(r)\sigma(\dd r)=\infty \quad\Longrightarrow\quad \text{a.s.} \quad \leb^1([0,\infty)\backslash L_{\ph})=0.
\end{split}\end{equation*}
\end{lem}

\begin{proof}
Assume that $\int_{0^{+}} \ph(r)\sigma(\dd r)<\infty$ and let $n\in\N$ and $\eps\in (0,1]$. The compensation formula for Poisson point processes yields
\[
\esp\left[\int_{0\leq t\leq n \atop \|x\|\leq\eps} \ph(\|x\|) J(\dd t,\dd x)\right]=n\int_{\|x\|\leq\eps} \ph(\|x\|) \pi(\dd x)=n\int_{(0,\eps]} \ph(r)\sigma(\dd r),
\]
see~\cite[p.~7]{Bertoin:1996uq}. Since the last integral tends to zero as $\eps\to 0$, the result follows from Fatou's lemma.

Conversely, suppose that $\int_{0^{+}} \ph(r)\sigma(\dd r)=\infty$ and let $n\in\N$ and $t_{0}\in [0,n]$. In addition, assume that $t_{0}\not\in L_{\ph}$. Then there are finitely many $s\in S_{1}$ with $|t_{0}-s|<\ph(\|\Delta X_{s}\|)$. In particular, there is a positive integer $m$ such that for every integer $m'\geq m+1$, the Borel set
\[
B_{m,m'}=\left\{ (t,x)\in\HH \:\biggl|\: |t-t_{0}|<\ph(\|x\|) \quad\text{and}\quad \frac{1}{m'}<\|x\|<\frac{1}{m} \right\}
\]
contains no pair $(t,\Delta X_{t})$ with $t\in (0,\infty)$. Hence $J(B_{m,m'})=0$. As $J$ is a Poisson measure with intensity $\leb^1\restr{(0,\infty)}\otimes\pi$, this can happen with probability at most
\[
\exp\left(-\int_{\frac{1}{m'}<\|x\|<\frac{1}{m}} \ph(\|x\|) \pi(\dd x)\right)=\exp\left(-\int_{\left(\frac{1}{m'},\frac{1}{m}\right)} \ph(r) \sigma(\dd r)\right),
\]
which tends to zero as $m'\to\infty$. As a result, any real number $t_{0}\in [0,n]$ belongs to $L_{\ph}$ with probability one. Fubini's theorem leads to
\[
n=\int_{0}^n \prob( t_{0}\in L_{\ph} ) \dd t_{0}=\esp\left[\leb^1(L_{\ph}\cap [0,n])\right]
\]
and the result follows.
\end{proof}

The following lemma will be called upon at various points of the proof of Theorem~\ref{grintisoholdlevy}. We omit its proof because it is a straightforward consequence of the fact that $J$ is a Poisson measure with intensity $\leb^1\restr{(0,\infty)}\otimes\pi$.

\begin{lem}\label{nbrsautsfini}
With probability one, for every finite interval $I$ and every positive real $\eps$, there are finitely many jump times $s\in S_{1}\cap I$ such that $\|\Delta X_{s}\|>\eps$.
\end{lem}

The proof of Theorem~\ref{grintisoholdlevy} uses some techniques developed in~\cite{Durand:2007uq} which we now recall in general dimension $d$. Let $V$ denote a nonempty open subset of $\R^d$. Let $I$ be a denumerable (i.e. countably infinite) set and let $(x_{i},r_{i})_{i\in I}$ denote a family in $\R^d\times (0,\infty)$ enjoying
\[
\sup_{i\in I}r_{i}<\infty \qquad\text{and}\qquad \forall m\in\N \qquad \#\left\{ i\in I \:\biggl|\: \|x_{i}\|<m \text{ and } r_{i}>\frac{1}{m} \right\}<\infty.
\]
The family $(x_{i},r_{i})_{i\in I}$ is called a {\em homogeneous ubiquitous system} in $V$ if for Lebesgue almost every point $x\in V$, there are infinitely many $i\in I$ such that $\|x-x_{i}\|<r_{i}$. Moreover, for every gauge function $g\in\jauge_{d}$, let $(g^{1/d})^{-1}$ denote the pseudo-inverse of $g^{1/d}$, which is defined by
\[
(g^{1/d})^{-1}(r)=\inf\{ \rho\in [0,\eps_{g}) \:|\: g^{1/d}(\rho)\geq r \}
\]
for every $r\in [0,\sup_{[0,\eps_{g})} g^{1/d})$. The following result is proven in~\cite{Durand:2007uq}.

\begin{thm}\label{UBIQUITYHOM}
Let $V$ denote a nonempty open subset of $\R^d$, let $(x_{i},r_{i})_{i\in I}$ denote a homogeneous ubiquitous system in $V$ and let $g\in\jauge_{d}$. Then, for every nonnegative nondecreasing function $\psi:[0,\infty)\to\R$ that coincides with $(g^{1/d})^{-1}$ in a neighborhood of the origin,
\[
\left\{ x\in\R^d \:\bigl|\: \|x-x_{i}\|<\psi(r_{i})\text{ for infinitely many }i\in I \right\}\in\grint^g(V).
\]
\end{thm}

\begin{rem}
For every continuous nondecreasing function $\ph:[0,\infty)\to\R$ enjoying $\ph(0)=0$ and $\int_{0^{+}} \ph(r)\sigma(\dd r)=\infty$, the family $(s,\ph(\|\Delta X_{s}\|))_{s\in S_{1}}$ is almost surely a homogeneous ubiquitous system in $(0,\infty)$. Indeed, Lemma~\ref{sizelph} ensures that the set $L_{\ph}$ almost surely has full Lebesgue measure in $(0,\infty)$. In addition, the set $S_{1}$ is almost surely denumerable. Furthermore, $\ph(\|\Delta X_{s}\|)\leq \ph(1)$ for all $s\in S_{1}$ and, for every $m\in\N$, the set of all times $s\in S_{1}$ such that $|s|<m$ and $\ph(\|\Delta X_{s}\|)>1/m$ is almost surely finite owing to Lemma~\ref{nbrsautsfini}. Theorem~\ref{UBIQUITYHOM} can therefore be applied to the family $(s,\ph(\|\Delta X_{s}\|))_{s\in S_{1}}$.
\end{rem}

\subsection{Proof of Theorem~\ref{grintisoholdlevy}}

For the sake of clarity, we split the statement of Theorem~\ref{grintisoholdlevy} into five propositions, namely, Propositions~\ref{prphauuveh0hau} to~\ref{hauuvehinfty}, which we present and prove all along this section.

Before stating these propositions, we begin by fixing some notations and making some remarks. Assume that~(\ref{condcjlevy}) holds and let $u\in\jauge_{1}$ with $h_{u}<\infty$. In addition, let $\tilde u$ denote a continuous nondecreasing function defined on $[0,\infty)$ that coincides with $u$ in a neighborhood of zero. As $h_{u}=h_{\tilde u}$, Lemma~\ref{sizelph} ensures that for every fixed $h\in (0,h_{u})$, with probability one, for every $n\in\N$,
\[
A_{n,\eps}(h)=\int_{0\leq t\leq n \atop \|x\|\leq\eps} \tilde u(\|x\|^{1/h}) J(\dd t,\dd x)
\]
tends to zero as $\eps\to 0$. Hence with probability one, for every $m\in\N$ such that $h_{u}-1/m>0$ and every $n\in\N$, the integral $A_{n,\eps}(h_{u}-1/m)$ tends to zero as $\eps\to 0$. Moreover, for every $h\in (0,h_{u})$, there is a positive integer $m$ enjoying $h\leq h_{u}-1/m$. Thus $A_{n,\eps}(h)\leq A_{n,\eps}(h_{u}-1/m)$ because $h\mapsto A_{n,\eps}(h)$ is nondecreasing. As a result,
\begin{equation}\label{sizelphconv}
\forall h\in (0,h_{u}) \quad \forall n\in\N \qquad \lim_{\eps\downarrow 0}\downarrow \int_{0\leq t\leq n \atop \|x\|\leq\eps} \tilde u(\|x\|^{1/h}) J(\dd t,\dd x)=0.
\end{equation}
Likewise, one can establish thanks to Lemma~\ref{sizelph} that with probability one,
\begin{equation}\label{sizelphdiv}
\forall h\in (h_{u},\infty) \qquad \leb^1([0,\infty)\backslash L_{\tilde u\circ\Id^{1/h}})=0.
\end{equation}
It follows that~(\ref{sizelphconv}),~(\ref{sizelphdiv}) and the statements of Proposition~\ref{locholdlevy} and Lemma~\ref{nbrsautsfini} simultaneously hold with probability one. From now on, we assume that the corresponding event occurs.

Consider $v\in\jaugew_{1}$ and let $\tilde v$ be a continuous nondecreasing function defined on $[0,\infty)$ that coincides with $v$ in a neighborhood of zero. The limit
\[
\gamma_{v}=\lim_{r\to 0}\frac{\log v(r)}{\log r}
\]
exists and belongs to the interval $[0,1]$. Furthermore, some routine calculations show that $h_{u\circ v}=h_{u}\gamma_{v}$.

Let $h\in [0,1/\beta')$ and let $V$ be a nonempty open subset of $(0,\infty)$.

\begin{prp}\label{prphauuveh0hau}
If $h<h_{u\circ v}$, then $\hau^{u\circ v}(E_{h}\cap V)=\hau^{u\circ v}(\tilde E_{h}\cap V)=0$.
\end{prp}

\begin{proof}
Let $\alpha\in (h,h_{u\circ v})$. Proposition~\ref{locholdlevy} implies that $E_{h}\backslash S\subseteq\tilde E_{h}\subseteq L_{\Id^{1/\alpha}}$. As $S$ is countable, $\hau^{u\circ v}(E_{h})\leq\hau^{u\circ v}(L_{\Id^{1/\alpha}})$ and $\hau^{u\circ v}(\tilde E_{h})\leq\hau^{u\circ v}(L_{\Id^{1/\alpha}})$. Hence, it suffices to show that the set $L_{\Id^{1/\alpha}}$ has zero Hausdorff $u\circ v$-measure.

To this end, observe that $\alpha<h_{u}\gamma_{v}$ so that $\alpha/(\gamma_{v}-\eta)<h_{u}$ for some $\eta\in (0,\gamma_{v})$. In addition, let $n\in\N$ and let $\eps$ denote a positive real small enough to ensure that $u\circ v(r)=\tilde u\circ\tilde v(r)\leq\tilde u(r^{\gamma_{v}-\eta})$ for all $r\in [0,\eps]$. For $t\in L_{\Id^{1/\alpha}}\cap [0,n]$, there are infinitely many jump times $s\in S_{1}$ such that $|t-s|<\|\Delta X_{s}\|^{1/\alpha}$. Such a jump time $s$ necessarily belongs to $[0,n+1]$ since $\|\Delta X_{s}\|\leq 1$. Furthermore, there are finitely many $s\in S_{1}\cap [0,n+1]$ enjoying $\|\Delta X_{s}\|>\eps^\alpha$ because of Lemma~\ref{nbrsautsfini}. As a consequence, there exists a jump time $s\in S_{1}\cap [0,n+1]$ such that $\|\Delta X_{s}\|\leq\eps^\alpha$ and $|t-s|<\|\Delta X_{s}\|^{1/\alpha}$. Hence
\[
L_{\Id^{1/\alpha}}\cap [0,n]\subseteq \bigcup_{s\in S_{1}\cap [0,n+1] \atop \|\Delta X_{s}\|\leq\eps^\alpha} (s-\|\Delta X_{s}\|^{1/\alpha},s+\|\Delta X_{s}\|^{1/\alpha}).
\]
This covering yields
\begin{equation*}\begin{split}
\hau^{u\circ v}_{\eps}(L_{\Id^{1/\alpha}}\cap [0,n]) &\leq 2\sum_{s\in S_{1}\cap [0,n+1] \atop \|\Delta X_{s}\|\leq\eps^\alpha} u\circ v(\|\Delta X_{s}\|^{1/\alpha}) \\
&\leq 2\int_{0\leq t\leq n+1 \atop \|x\|\leq\eps^\alpha} \tilde u\bigl(\|x\|^{\frac{\gamma_{v}-\eta}{\alpha}}\bigr) J(\dd t,\dd x).
\end{split}\end{equation*}
Owing to~(\ref{sizelphconv}), this integral tends to zero as $\eps\to 0$. Proposition~\ref{prphauuveh0hau} follows.
\end{proof}

\begin{prp}\label{prphauuveh0grint}
If $h<h_{u\circ v}$, then $\tilde E_{h}\not\in\grint^{u\circ v}(V)$.
\end{prp}

\begin{proof}
Let us build a gauge function $\overline{u}$ in $\jauge_{1}$ such that $\overline{u}\prec u$ and $h_{\overline{u}}\geq h_{u}$. For all $n\in\N$, let $\alpha_{n}=(1+1/n)/h_{u}$. Note that the sequence $(\alpha_{n})_{n\in\N}$ decreases and converges to $1/h_{u}$. Moreover, $\lim_{0^+}u=0$ and $\int_{0^+} u(r^{\alpha_{n}})\sigma(\dd r)<\infty$ for all $n\in\N$. So there exists a decreasing sequence $(r_{n})_{n\in\N}$ in $(0,1]$ such that the functions $u$ and $r\mapsto u(r)/r$ are respectively nondecreasing on $[0,r_{1}]$ and nonincreasing on $(0,r_{1}]$ and such that
\[
u(r_{n})\leq u(r_{n-1})e^{-1/n} \qquad\text{and}\qquad \int_{(0,{r_{n-1}}^{1/\alpha_{1}}]} u(r^{\alpha_{n}})\sigma(\dd r)\leq\frac{1}{(n+1)^3}
\]
for every integer $n\geq 2$. Observe that $(r_{n})_{n\in\N}$ necessarily converges to zero since $u(r_{n})$ tends to zero as $n\to\infty$ and $u$ is positive and continuous on $(0,r_{1}]$. For every integer $n\geq 2$ and every $r\in (r_{n},r_{n-1}]$, let
\[
w(r)=n+\frac{\log u(r_{n-1})-\log u(r)}{\log u(r_{n-1})-\log u(r_{n})}.
\]
The function $w$ is continuous and nonincreasing on $(0,r_{1}]$. For every $r$ in this interval, let $\overline{u}(r)=u(r)w(r)$. Consider an integer $n\geq 2$ and two reals $r$ and $r'$ such that $r_{n}<r\leq r'\leq r_{n-1}$. Then $\overline{u}(r')-\overline{u}(r)$ vanishes if $u(r)=u(r')$ and otherwise it is at least
\[
\left(u(r')-u(r)\right)n\left(1-\frac{\log\frac{u(r')}{u(r)}}{\frac{u(r')}{u(r)}-1}\cdot\frac{1}{n\log\frac{u(r_{n-1})}{u(r_{n})}}\right)\geq 0.
\]
The function $\overline{u}$ is therefore nondecreasing on $(r_{n},r_{n-1}]$ for each integer $n\geq 2$. As it is continuous on $(0,r_{1}]$, it is nondecreasing on this interval. Given $n_{0}\in\N$, observe that $u(r^{\alpha_{n_{0}}})\leq u(r^{\alpha_{n}})$ and $w(r^{\alpha_{n_{0}}})\leq n+1$ for every integer $n\geq n_{0}+1$ and every real $r\in ({r_{n}}^{1/\alpha_{n_{0}}},{r_{n-1}}^{1/\alpha_{n_{0}}}]$, so that
\begin{equation*}\begin{split}
\int_{(0,{r_{n_{0}}}^{1/\alpha_{n_{0}}}]} \overline{u}(r^{\alpha_{n_{0}}})\sigma(\dd r)
&=\sum_{n=n_{0}+1}^\infty\int_{({r_{n}}^{1/\alpha_{n_{0}}},{r_{n-1}}^{1/\alpha_{n_{0}}}]} u(r^{\alpha_{n_{0}}})w(r^{\alpha_{n_{0}}})\sigma(\dd r) \\
&\leq\sum_{n=n_{0}+1}^\infty (n+1)\int_{(0,{r_{n-1}}^{1/\alpha_{1}}]} u(r^{\alpha_{n}})\sigma(\dd r) \\
&\leq\sum_{n=n_{0}+1}^\infty\frac{1}{(n+1)^2}<\infty.
\end{split}\end{equation*}
This implies that $\overline{u}$ tends to zero at zero because $\sigma$ has infinite total mass, owing to the finiteness of $h_{u}$. Hence $\overline{u}\in\jauge_{1}$ and $\overline{u}\prec u$. This also implies that $1/\alpha_{n_{0}}\leq h_{\overline{u}}$. Letting $n_{0}\to\infty$ yields $h_{u}\leq h_{\overline{u}}$. Besides, $h_{\overline{u}}\leq h_{u}<\infty$ because $\overline{u}\prec u$. Then the analog of~(\ref{sizelphconv}) with $\overline{u}$ instead of $u$ holds with probability one. Assume that the corresponding event occurs. As $h<h_{u\circ v}=h_{u}\gamma_{v}\leq h_{\overline{u}}\gamma_{v}=h_{\overline{u}\circ v}$, the conclusion of Proposition~\ref{prphauuveh0hau} still holds with $\overline{u}$ instead of $u$. Thus $\hau^{\overline{u}\circ v}(\tilde E_{h}\cap V)=0$. Proposition~\ref{prphauuveh0grint} then follows from Theorem~\ref{GRINTSTABLE} and the fact that $\overline{u}\circ v\prec u\circ v$.
\end{proof}

\begin{prp}\label{hauuvtildeehinftygrint}
If $h\geq h_{u\circ v}$, then $\tilde E_{h}\in\grint^{u\circ v}(V)$.
\end{prp}

\begin{proof}
Let $\alpha\in (h,1/\beta]$. Then $\alpha>h_{u\circ v}=h_{u}\gamma_{v}$ so that $\alpha/(\gamma_{v}+\eps)>h_{u}$ for every $\eps\in (0,-\gamma_{v}+\alpha/h_{u})$ with the convention that the upper bound of this interval is infinite if $h_{u}$ vanishes. For every positive $r$ small enough, we have $\tilde v(r)\geq r^{\gamma_{v}+\eps}$ and hence $\tilde u\circ\tilde v(r^{1/\alpha})\geq\tilde u(r^{(\gamma_{v}+\eps)/\alpha})$. Consequently, the set $L_{\tilde u\circ\tilde v\circ \Id^{1/\alpha}}$ contains $L_{\tilde u\circ\Id^{(\gamma_{v}+\eps)/\alpha}}$. Meanwhile, this last set has full Lebesgue measure in $[0,\infty)$ owing to~(\ref{sizelphdiv}). It follows that $L_{\tilde u\circ\tilde v\circ\Id^{1/\alpha}}$ has full Lebesgue measure in $(0,\infty)$.

Along with Lemma~\ref{nbrsautsfini}, this result ensures that the family $(s,\tilde u\circ\tilde v(\|\Delta X_{s}\|^{1/\alpha}))_{s\in S_{1}}$ is a homogeneous ubiquitous system in $(0,\infty)$. Theorem~\ref{UBIQUITYHOM} implies that for any nonnegative nondecreasing function $\ph:[0,\infty)\to\R$ that coincides with $(\tilde u\circ\tilde v)^{-1}$ in a neighborhood of zero, the set of all reals $t$ such that $|t-s|<\ph(\tilde u\circ\tilde v(\|\Delta X_{s}\|^{1/\alpha}))$ for infinitely many $s\in S_{1}$ belongs to the class $\grint^{\tilde u\circ\tilde v}((0,\infty))$ hence to the class $\grint^{u\circ v}((0,\infty))$ as $\tilde u\circ\tilde v$ and $u\circ v$ coincide in a neighborhood of the origin. The $G_{\delta}$-set of all reals $t$ such that $|t-s|<\|\Delta X_{s}\|^{1/\alpha}$ for infinitely many $s\in S_{1}$ contains the aforementioned set because $\ph(\tilde u\circ\tilde v(r))\leq r$ for every positive $r$ small enough. Thus it belongs to the class $\grint^{u\circ v}((0,\infty))$ as well. The interval $[0,\infty)$ also belongs to this class. Theorem~\ref{GRINTSTABLE} finally ensures that
\[
\forall \alpha\in (h,1/\beta] \qquad L_{\Id^{1/\alpha}}\in\grint^{u\circ v}((0,\infty)).
\]

In addition, as $S$ is countable, $\R\backslash S$ is a $G_{\delta}$-set of full Lebesgue measure in $(0,\infty)$, thereby belonging to $\grint^{u\circ v}((0,\infty))$ thanks to Theorem~\ref{GRINTSTABLE}. Moreover, Proposition~\ref{locholdlevy} and the observation that $\alpha\mapsto L_{\Id^{1/\alpha}}$ is nondecreasing yield
\[
\tilde E_{h}=\left(\R\backslash S\right)\cap\bigcap_{n\in\N \atop h+1/n\leq 1/\beta} L_{\Id^{1/(h+1/n)}}
\]
so that $\tilde E_{h}$ is a countable intersection of sets of the class $\grint^{u\circ v}((0,\infty))$. Theorem~\ref{GRINTSTABLE} then implies that $\tilde E_{h}$ belongs to $\grint^{u\circ v}((0,\infty))$. Proposition~\ref{hauuvtildeehinftygrint} follows.
\end{proof}

\begin{prp}\label{hauuvtildeehinftyhau}
If $h\geq h_{u\circ v}$, then $\hau^{u\circ v}(\tilde E_{h}\cap V)=\infty$.
\end{prp}

\begin{proof}
We begin by assuming that $u\not\prec\Id$. One easily checks that $1/\beta\leq h_{u}<\infty$ so that $\beta$ does not vanish. Hence the identity function belongs to $\jauge_{1}$ and enjoys $h_{\Id}=1/\beta<\infty$. Then the analog of~(\ref{sizelphdiv}) with $\Id$ instead of $u$ holds with probability one and we may suppose that the corresponding event occurs. Furthermore, since $\gamma_{v}/\beta\leq h_{u}\gamma_{v}=h_{u\circ v}\leq h<1/\beta$, we have $\gamma_{v}<1$ so that $v\prec\Id$. In consequence, the function $\underline{v}:r\mapsto v(r)/\log(v(r)/r)$ belongs to $\jaugew_{1}$ and enjoys $v\prec\underline{v}$ and $\gamma_{\underline{v}}=\gamma_{v}$. As $h\geq h_{u\circ v}=h_{u}\gamma_{v}\geq h_{\Id}\gamma_{\underline{v}}=h_{\Id\circ\underline{v}}$, the conclusion of Proposition~\ref{hauuvtildeehinftygrint} still holds with $\Id$ instead of $u$ and $\underline{v}$ instead of $v$. Thus $\tilde E_{h}\in\grint^{\underline{v}}(V)$. Theorem~\ref{GRINTSTABLE} then ensures that $\hau^v(\tilde E_{h}\cap V)=\infty$. The result finally stems from the fact that $\hau^{u\circ v}\geq C\hau^v$ for some positive $C$ as $u\in\jauge_{1}$.

Let us now assume that $u\prec\Id$ and build a gauge function $\underline{u}\in\jauge_{1}$ such that $u\prec\underline{u}$ and $h_{\underline{u}}\leq h_{u}$. To this end, let $\alpha_{n}=1/(h_{u}+1/n)$ for every $n\in\N$ and observe that $\int_{0^+} u(r^{\alpha_{n}})\sigma(\dd r)=\infty$. Meanwhile, there exists $r_{1}\in (0,1]$ such that $u$ is continuous and nondecreasing on $[0,{r_{1}}^{\alpha_{1}}]$ and $\rho:r\mapsto u(r)/r$ is nonincreasing on $(0,{r_{1}}^{\alpha_{1}}]$. Moreover, $\rho$ tends to infinity at zero. It follows that for every integer $n\geq 2$, there is a real number $r_{n}\in (0,r_{n-1})$ such that
\[
\rho({r_{n}}^{\alpha_{n}})\geq\rho({r_{n-1}}^{\alpha_{n-1}})e^{1/n} \qquad\text{and}\qquad \int_{(r_{n},r_{n-1}]} u(r^{\alpha_{n}}) \sigma(\dd r)\geq 1.
\]
Note that the sequence $({r_{n}}^{\alpha_{n}})_{n\in\N}$ decreases and converges to zero. For every integer $n\geq 2$ and every real $r\in ({r_{n}}^{\alpha_{n}},{r_{n-1}}^{\alpha_{n-1}}]$, let
\[
w(r)=n+\frac{\log\rho(r)-\log\rho({r_{n-1}}^{\alpha_{n-1}})}{\log\rho({r_{n}}^{\alpha_{n}})-\log\rho({r_{n-1}}^{\alpha_{n-1}})}.
\]
Then, the function $w$ is continuous, nonincreasing and positive on $(0,{r_{1}}^{\alpha_{1}}]$. Let $\underline{u}(r)=u(r)/w(r)$ for all $r\in (0,{r_{1}}^{\alpha_{1}}]$. This function is continuous and nondecreasing on $(0,{r_{1}}^{\alpha_{1}}]$ and tends to zero at zero, so it belongs to $\jauge$. Furthermore, $u\prec\underline{u}$ since $w$ monotonically tends to infinity at zero. In addition, $\underline{u}$ belongs to $\jauge_{1}$. Indeed, $r\mapsto\underline{u}(r)/r$ is continuous at ${r_{n}}^{\alpha_{n}}$ for every integer $n\geq 2$ and is nonincreasing on $({r_{n}}^{\alpha_{n}},{r_{n-1}}^{\alpha_{n-1}}]$ because, for $r\leq r'$ in this interval, $\underline{u}(r)/r-\underline{u}(r')/r'$ vanishes if $\rho(r)=\rho(r')$ and is at least
\[
\frac{\rho(r)-\rho(r')}{w(r)w(r')}n\left(1-\frac{\log\frac{\rho(r)}{\rho(r')}}{\frac{\rho(r)}{\rho(r')}-1}\cdot\frac{1}{n\log\frac{\rho({r_{n}}^{\alpha_{n}})}{\rho({r_{n-1}}^{\alpha_{n-1}})}}\right)\geq 0
\]
otherwise. Let $n_{0}$ denote an integer greater than $1$. For every integer $n\geq n_{0}$ and every $r\in (r_{n},r_{n-1}]$, we have ${r_{n}}^{\alpha_{n}}<r^{\alpha_{n}}\leq r^{\alpha_{n_{0}}}$, so that $n+1=w({r_{n}}^{\alpha_{n}})\geq w(r^{\alpha_{n_{0}}})$ and $u(r^{\alpha_{n}})\leq u(r^{\alpha_{n_{0}}})$. Hence
\begin{equation*}\begin{split}
\int_{(0,r_{1}]}\underline{u}(r^{\alpha_{n_{0}}})\sigma(\dd r)&\geq\sum_{n=n_{0}}^\infty\int_{(r_{n},r_{n-1}]}\frac{u(r^{\alpha_{n_{0}}})}{w(r^{\alpha_{n_{0}}})}\sigma(\dd r)\\
&\geq\sum_{n=n_{0}}^\infty\frac{1}{n+1}\int_{(r_{n},r_{n-1}]}u(r^{\alpha_{n}})\sigma(\dd r)\geq \sum_{n=n_{0}}^\infty\frac{1}{n+1}=\infty.
\end{split}\end{equation*}
As a result, $h_{\underline{u}}\leq 1/\alpha_{n_{0}}$. Letting $n_{0}\to\infty$ yields $h_{\underline{u}}\leq h_{u}<\infty$. Thus the analog of~(\ref{sizelphdiv}) with $\underline{u}$ instead of $u$ holds with probability one and we may assume that the corresponding event occurs. Then, since $h\geq h_{u\circ v}=h_{u}\gamma_{v}\geq h_{\underline{u}}\gamma_{v}=h_{\underline{u}\circ v}$, the conclusion of Proposition~\ref{hauuvtildeehinftygrint} holds with $\underline{u}$ instead of $u$. Therefore, $\tilde E_{h}\in\grint^{\underline{u}\circ v}(V)$. As $u\circ v\prec \underline{u}\circ v$, Theorem~\ref{GRINTSTABLE} implies that $\hau^{u\circ v}(\tilde E_{h}\cap V)=\infty$.
\end{proof}

\begin{prp}\label{hauuvehinfty}
If $h\in [h_{u\circ v},h_{u}]$, then $\hau^{u\circ v}(E_{h}\cap V)=\infty$.
\end{prp}

\begin{proof}
We begin by supposing that $h=h_{u\circ v}$. In particular, $h\geq h_{u\circ v}$ and Proposition~\ref{hauuvtildeehinftyhau} ensures that $\hau^{u\circ v}(\tilde E_{h}\cap V)=\infty$. Moreover, Proposition~\ref{locholdlevy} along with the fact that $\alpha\mapsto L_{\Id^{1/\alpha}}$ is nondecreasing leads to
\[
E_{h}\backslash S=\tilde E_{h}\backslash\bigcup_{m\in\N \atop h-1/m>0} L_{\Id^{1/(h-1/m)}}.
\]
Let $m\in\N$ with $h-1/m>0$. Since $h-1/m<h=h_{u}\gamma_{v}$, there exists a real number $\eta\in (0,\gamma_{v})$ such that $(h-1/m)/(\gamma_{v}-\eta)<h_{u}$. In addition, let $n\in\N$ and let $\eps$ denote a positive real small enough to ensure that $u\circ v(r)\leq\tilde u(r^{\gamma_{v}-\eta})$ for all $r\in [0,\eps]$. The set $L_{\Id^{1/(h-1/m)}}\cap [0,n]$ is covered by the open intervals of center $s\in S_{1}\cap [0,n+1]$ and radius $\|\Delta X_{s}\|^{1/(h-1/m)}\leq\eps$. This covering leads to
\[
\hau^{u\circ v}_{\eps}(L_{\Id^{1/(h-1/m)}}\cap [0,n])\leq 2\int_{0\leq t\leq n+1 \atop \|x\|\leq\eps^{h-1/m}} \tilde u\bigl(\|x\|^{\frac{\gamma_{v}-\eta}{h-1/m}}\bigr) J(\dd t,\dd x).
\]
Owing to~(\ref{sizelphconv}), this integral tends to zero as $\eps\to 0$. It follows that $L_{\Id^{1/(h-1/m)}}$ has zero Hausdorff $u\circ v$-measure. Hence, Proposition~\ref{hauuvehinfty} holds for $h=h_{u\circ v}$.

Let us assume that $h>h_{u\circ v}$. Then, $h_{u}$ is necessarily positive. Let $w=\Id^{h/h_{u}}$. As $h\leq h_{u}$, we have $w\in\jaugew_{1}$ and $\gamma_{w}=h/h_{u}$ so that $h_{u\circ w}=h_{u}\gamma_{w}=h$. Using $w$ instead of $v$ in the first part of the proof, we obtain $\hau^{u\circ w}(E_{h}\cap V)=\infty$. Meanwhile, $\gamma_{v}<h/h_{u}=\gamma_{w}$ so $u\circ v(r)\geq u\circ w(r)$ for every positive real $r$ small enough. In consequence, $\hau^{u\circ v}\geq\hau^{u\circ w}$. Hence, Proposition~\ref{hauuvehinfty} holds for $h>h_{u\circ v}$.
\end{proof}

\subsection{Proof of Theorem~\ref{grintisoholdlevy1}}

Let us suppose that~(\ref{condcjlevy}) holds, let $g\in\jauge$ and $\tilde g$ denote a continuous nondecreasing function defined on $[0,\infty)$ that coincides with $g_{1}\in\jauge_{1}\cup\{0\}$ in a neighborhood of the origin.

To begin with, let us assume that $h_{g}=\infty$. By following the proof of~(\ref{sizelphconv}), it is easy to establish that with probability one,
\[
\forall h\in (0,\infty) \quad \forall n\in\N \qquad \lim_{\eps\downarrow 0}\downarrow \int_{0\leq t\leq n \atop \|x\|\leq\eps} \tilde g(\|x\|^{1/h}) J(\dd t,\dd x)=0.
\]
Let us suppose that the event on which this assertion and the statements of Proposition~\ref{locholdlevy} and Lemma~\ref{nbrsautsfini} hold occurs. Let $h\in [0,1/\beta')$, $\alpha\in (h,1/\beta]$, $n\in\N$ and let $\eps$ denote a positive real small enough to ensure that $g_{1}$ and $\tilde g$ coincide on $[0,\eps]$. Owing to Lemma~\ref{nbrsautsfini}, the set $L_{\Id^{1/\alpha}}\cap [0,n]$ is covered by the open intervals of center $s\in S_{1}\cap [0,n+1]$ and radius $\|\Delta X_{s}\|^{1/\alpha}\leq\eps$. This covering leads to
\[
\hau^{g_{1}}_{\eps}(L_{\Id^{1/\alpha}}\cap [0,n])\leq 2\int_{0\leq t\leq n+1 \atop \|x\|\leq\eps^\alpha} \tilde g(\|x\|^{1/\alpha}) J(\dd t,\dd x).
\]
Letting $\eps\to 0$ implies that the set $L_{\Id^{1/\alpha}}$ has zero Hausdorff $g_{1}$-measure. It follows from~(\ref{prpjauge1}) that $\hau^g(L_{\Id^{1/\alpha}})=0$. As $S$ is countable, Proposition~\ref{locholdlevy} shows that $E_{h}$ and $\tilde E_{h}$ have zero $g$-measure as well.

Conversely, let us assume that $h_{g}<\infty$. The function $g_{1}$ then necessarily belongs to $\jauge_{1}$. Let us apply Theorem~\ref{grintisoholdlevy} with $u=g_{1}$ and $v=\Id$. Then, with probability one, for every real $h\in [0,1/\beta')$ and every nonempty open subset $V$ of $(0,\infty)$,
\[
\hau^{g_{1}}(E_{h}\cap V)=\begin{cases}
0 & \text{if }h<h_{g_{1}} \\
\infty & \text{if }h=h_{g_{1}}
\end{cases}
\qquad\text{and}\qquad
\hau^{g_{1}}(\tilde E_{h}\cap V)=\begin{cases}
0 & \text{if }h<h_{g_{1}} \\
\infty & \text{if }h\geq h_{g_{1}}
\end{cases}.
\]
Theorem~\ref{grintisoholdlevy1} now follows from~(\ref{prpjauge1}) along with the observation that $h_{g}=h_{g_{1}}$.

\subsection{Proof of Theorem~\ref{grintisoholdlevy2}}

Let us suppose that $\beta>0$ and that~(\ref{condcjlevy}) holds. As $h_{\Id}=1/\beta<\infty$, we can apply Theorem~\ref{grintisoholdlevy} with $u=\Id$. Hence, with probability one, for every $g\in\jaugew_{1}$, every $h\in [0,1/\beta')$ and every nonempty open $V\subseteq (0,\infty)$,
\[
\hau^g(E_{h}\cap V)=\hau^g(\tilde E_{h}\cap V)=\begin{cases}
0 & \text{si }h<h_{g} \\
\infty & \text{si }h\geq h_{g}
\end{cases}.
\]
Theorem~\ref{grintisoholdlevy2} is then a direct consequence of~(\ref{prpjauge1}) and the fact that $h_{g}=h_{g_{1}}$.

\subsection{Proof of Theorem~\ref{grintisoholdlevy3}}

Assume that~(\ref{condcjlevy}) holds. In order to prove the first part of Theorem~\ref{grintisoholdlevy3}, let $g\in\jauge_{1}$. We begin by supposing that $h_{g}=\infty$. Hence, there exists a gauge function $\overline{g}\in\jauge_{1}$ such that $\overline{g}\prec g$ and $h_{\overline{g}}=\infty$. Applying Theorem~\ref{grintisoholdlevy1} with $\overline{g}$ rather than $g$, one proves that with probability one, for each $h\in [0,1/\beta')$, the set $\tilde E_{h}$ has zero Hausdorff $\overline{g}$-measure. The result then follows from Theorem~\ref{GRINTSTABLE}. Conversely, if $h_{g}<\infty$, the result is easily obtained by applying Theorem~\ref{grintisoholdlevy} with $u=g$ and $v=\Id$.

If $\beta>0$, then $h_{\Id}=1/\beta<\infty$. So the second part of Theorem~\ref{grintisoholdlevy3} is straightforwardly obtained by applying Theorem~\ref{grintisoholdlevy} with $u=\Id$.

\section{Moduli of continuity}\label{grintlevymodc}

This section is devoted to the proof of Theorem~\ref{grintmodclevy}. Let $g\in\jauge$ and $w\in\modc$ with $\int_{0^{+}} g_{1}(w^{-1}(r))\sigma(\dd r)=\infty$. Note that $g_{1}$ belongs to $\jauge_{1}$. Let $\tilde g$ denote a continuous nondecreasing function defined on $[0,\infty)$ that coincides with $g_{1}$ in a neighborhood of zero and let $\tilde w$ be a continuous increasing function defined on $[0,\infty)$ that tends to infinity at infinity and coincides with $w$ in a neighborhood of the origin.

Observe that there are two real numbers $\kappa_{1}$ and $\kappa_{2}$ with $1<\kappa_{1}\leq\kappa_{2}$ and a positive real $\delta_{0}$ such that $\kappa_{1}\tilde w(\delta)\leq\tilde w(2\delta)\leq\kappa_{2}\tilde w(\delta)$ for all $\delta\in [0,\delta_{0}]$. Let $\ph_{q}:r\mapsto \tilde w^{-1}(r/{\kappa_{1}}^q)$ for each $q\in\N$ and consider
\[
\tilde F_{w}=\left(\bigcap_{q=1}^\infty\downarrow L_{\ph_{q}}\right)\backslash S.
\]
Theorem~\ref{grintmodclevy} is then a straightforward consequence of the two following lemmas.

\begin{lem}
With probability one, $\tilde F_{w}\in\grint^{g_{1}}((0,\infty))$ and $\hau^g(\tilde F_{w}\cap V)=\hau^g(V)$ for every open subset $V$ of $(0,\infty)$.
\end{lem}

\begin{proof}
As $\kappa_{1}\tilde w(\delta)\leq\tilde w(2\delta)$ for any $\delta\in [0,\delta_{0}]$, we have $\tilde w^{-1}(r/\kappa_{1})\geq\tilde w^{-1}(r)/2$ for every $r>0$ small enough. Thus, for $q\in\N$ and $r$ small enough, $\ph_{q}(r)\geq\tilde w^{-1}(r)/2^q$, so $\tilde g(\ph_{q}(r))\geq \tilde g(\tilde w^{-1}(r))/2^q$ owing to the fact that $\tilde g\in\jauge_{1}$. In consequence, for every $\eps>0$ small enough,
\[
\int_{(0,\eps]} \tilde g(\ph_{q}(r))\sigma(\dd r)\geq \frac{1}{2^q} \int_{(0,\eps]} \tilde g(\tilde w^{-1}(r)) \sigma(\dd r)=\infty.
\]
Lemmas~\ref{sizelph} and~\ref{nbrsautsfini} then show that the family $(s,\tilde g\circ\ph_{q}(\|\Delta X_{s}\|))_{s\in S_{1}}$ is a homogeneous ubiquitous system in $(0,\infty)$ with probability one. As the functions $g_{1}$ and $\tilde g$ coincide in a neighborhood of the origin, Theorem~\ref{UBIQUITYHOM} ensures that $L_{\ph_{q}}$ almost surely belongs to the class $\grint^{g_{1}}((0,\infty))$. Furthermore, the set $S$ is almost surely countable by Lemma~\ref{nbrsautsfini}, so $\R\backslash S$ almost surely belongs to this class by Theorem~\ref{GRINTSTABLE}. This theorem finally implies that $\tilde F_{w}$ almost surely belongs to $\grint^{g_{1}}((0,\infty))$ as well.

Let us suppose that $g_{1}\prec\Id$. There exists a function $\underline{g}\in\jauge_{1}$ such that $g_{1}\prec\underline{g}$ and $\int_{0^{+}} \underline{g}(w^{-1}(r))\sigma(\dd r)=\infty$. By using $\underline{g}$ rather than $g_{1}$ above, it is easy to check that $\tilde F_{w}$ almost surely belongs to $\grint^{\underline{g}}((0,\infty))$. Theorem~\ref{GRINTSTABLE} then ensures that with probability one, $\hau^{g_{1}}(\tilde F_{w}\cap V)=\infty=\hau^{g_{1}}(V)$ for every open subset $V$ of $(0,\infty)$ and~(\ref{prpjauge1}) leads to $\hau^g(\tilde F_{w}\cap V)=\hau^g(V)$.

Conversely, let us assume that $g_{1}\not\prec\Id$. Hence, $\tilde g(r)=g_{1}(r)\leq C r$ for some $C>0$ and every $r>0$ small enough. Therefore, $\hau^{g_{1}}\leq C\hau^1$ and, by virtue of~(\ref{prpjauge1}), $\hau^g$ is a translation invariant Borel measure which is finite on compacts. It follows that $\hau^g=\eta\leb^1$ on the Borel subsets of $\R$ for some $\eta>0$. Meanwhile, for $q\in\N$ and $\eps>0$ small enough,
\[
\int_{(0,\eps]} \ph_{q}(r)\sigma(\dd r)\geq\frac{1}{C}\int_{(0,\eps]} \tilde g(\ph_{q}(r))\sigma(\dd r)=\infty
\]
so that the set $L_{\ph_{q}}$ almost surely has full Lebesgue measure in $(0,\infty)$ by Lemma~\ref{sizelph}. In addition, $S$ is almost surely countable by Lemma~\ref{nbrsautsfini}. Thus $\tilde F_{w}$ almost surely has full Lebesgue measure in $(0,\infty)$. As a result, with probability one, for every open subset $V$ of $(0,\infty)$, we have $\hau^g(\tilde F_{w}\cap V)=\eta\leb^1(\tilde F_{w}\cap V)=\eta\leb^1(V)=\hau^g(V)$.
\end{proof}

\begin{lem}
With probability one, $\tilde F_{w}\subseteq F_{w}$.
\end{lem}

\begin{proof}
Let us assume that the event on which the statement of Lemma~\ref{nbrsautsfini} holds occurs and let $t_{0}\in\tilde F_{w}$. Moreover, let us suppose that $t_{0}$ does not belong to $F_{w}$. Hence $w$ is a modulus of continuity of $X$ at $t_{0}$ so there are two positive reals $c$ and $\delta$ and a $d$-tuple $P$ of polynomials such that for every $t\geq 0$,
\[
|t - t_{0}| \leq \delta \qquad\Longrightarrow\qquad \| X_{t} - P(t - t_{0}) \| \leq c\ w(| t - t_{0} |).
\]
We may actually assume that for every $t\geq 0$,
\[
|t - t_{0}| \leq 1 \qquad\Longrightarrow\qquad \| X_{t} - P(t - t_{0}) \| \leq c\ \tilde w(| t - t_{0} |).
\]
Consider an integer $q>\log(3c\kappa_{2})/\log\kappa_{1}$. As $t_{0}\in L_{\ph_{q}}\backslash S$, there exists an injective sequence $(s_{n})_{n\in\N}$ in $S_{1}$ such that $0<|t_{0}-s_{n}|<\ph_{q}(\|\Delta X_{s_{n}}\|)$ for all $n\in\N$. Owing to Lemma~\ref{nbrsautsfini}, there are finitely many integers $n$ enjoying $\ph_{q}(\|\Delta X_{s_{n}}\|)>(1/2)\wedge\delta_{0}$. Thus $|t_{0}-s_{n}|<(1/2)\wedge\delta_{0}$ for some $n\in\N$. Let us suppose that $\|X_{t}-P(t-t_{0})\|<\|\Delta X_{s_{n}}\|/3$ for every $t\geq 0$ with $|t-t_{0}|\leq 2|s_{n}-t_{0}|$. For every integer $p$ large enough, we have $|s_{n}-1/p-t_{0}|\leq 2|s_{n}-t_{0}|$ so that
\begin{equation*}\begin{split}
&\left\|(X_{s_{n}}-P(s_{n}-t_{0}))-\left(X_{s_{n}-\frac{1}{p}}-P\left(s_{n}-\frac{1}{p}-t_{0}\right)\right)\right\| \\
\leq&\left\|X_{s_{n}}-P(s_{n}-t_{0})\right\|+\left\|X_{s_{n}-\frac{1}{p}}-P\left(s_{n}-\frac{1}{p}-t_{0}\right)\right\|<\frac{2}{3}\|\Delta X_{s_{n}}\|.
\end{split}\end{equation*}
Meanwhile, the left-hand side tends to $\|\Delta X_{s_{n}}\|$ as $p\to\infty$, so we end up with a contradiction. Hence there exists a nonnegative real $t$ satisfying $|t-t_{0}|\leq 2|s_{n}-t_{0}|$ and $\|X_{t}-P(t-t_{0})\|\geq\|\Delta X_{s_{n}}\|/3$. It follows that
\[
\frac{\|\Delta X_{s_{n}}\|}{3}\leq\|X_{t}-P(t-t_{0})\|\leq c\ \tilde w(|t-t_{0}|)\leq c\ \tilde w(2|s_{n}-t_{0}|)\leq c\kappa_{2} \ \tilde w(|s_{n}-t_{0}|).
\]
Since $\|\Delta X_{s_{n}}\|>{\kappa_{1}}^q\tilde w(|t_{0}-s_{n}|)$, we obtain ${\kappa_{1}}^q<3c\kappa_{2}$, which is a contradiction. As a result, $t_{0}\in F_{w}$. The set $\tilde F_{w}$ is thus included in $F_{w}$.
\end{proof}

\section{Lacunary wavelet series}\label{grintlevylac}

The methods developed in the previous sections enable us to investigate the size and large intersection properties of the singularity sets of a model of random wavelet series which generalizes that studied by S.~Jaffard in~\cite{Jaffard:2000mk}. The relevance of this model is due to the fact that many signals, images or mathematical functions can be represented in a wavelet basis using very few nonvanishing coefficients. Examples include the piecewise smooth functions, the images denoised via wavelet thresholding and the solutions of certain nonlinear hyperbolic equations, see~\cite{DeVore:1990rt,Donoho:1995uq,Donoho:1995kx}. The process that we consider is a wavelet series with only a given number of nonvanishing coefficients at each scale $j$. These coefficients are equal to $2^{-\underline{h}j}$ for some fixed $\underline{h}>0$ and their positions are chosen uniformly and independently.

\subsection{Presentation of the model}

Let $\N_{0}$ denote the set of all nonnegative integers. Let $\I=\intn{1}{2^d-1}$ and let $\psi^i$, $i\in\I$, denote wavelets in the Schwartz class such that the functions $x\mapsto 2^{dj/2}\psi^i(2^jx-k)$, for $(i,j,k)\in\I\times\Z\times\Z^d$, form an orthonormal basis of $L^2(\R^d)$, see~\cite{Lemarie-Rieusset:1986lr}. Since we are interested in pointwise properties, it is more convenient to work with wavelets on the $d$-dimensional torus $\T^d=\R^d/\Z^d$. Let $\phi:\R^d\to\T^d$ be the canonical surjection and let $\dist$ denote the quotient distance on $\T^d$. Let $\Lambda$ denote the set of all dyadic cubes of the torus, that is, sets of the form $\lambda=\phi(2^{-j}(k+[0,1)^d))$ with $j\in\N_{0}$ and $k\in\intn{0}{2^j-1}^d$. Let $\gene{\lambda}=j$ be the generation of $\lambda$ et $x_{\lambda}=\phi(k 2^{-j})$. For all $i\in\I$ and $\lambda=\phi(2^{-j}(k+[0,1[^d))\in\Lambda$, let $\Psi^i_{\lambda}$ denote the function in $L^2(\T^d)$ that corresponds to the $\Z^d$-periodic function
\[
x\mapsto\sum_{m\in\Z^d} \psi^i\left(2^j(x-m)-k\right).
\]
The functions $2^{d\gene{\lambda}/2} \Psi^i_{\lambda}$, along with the constant function equal to one on $\T^d$, then form a wavelet basis of $L^2(\T^d)$. We refer to~\cite{Meyer:1992lr} for details.

Let $\underline{h}\in (0,\infty)$ and pick an integer $m_{i,j}\in\intn{0}{2^{dj}}$ for every $i\in\I$ and every $j\in\N_{0}$. The process that we consider is a wavelet series, denoted by $R$, with only $m_{i,j}$ nonvanishing coefficients at each scale $j$ and in each direction $i$, for $j$ large enough. These coefficients are equal to $2^{-\underline{h}j}$ and their positions are chosen uniformly and independently. To be specific, let $(X_{i,n})_{(i,n)\in\I\times\N_{0}}$ denote a family of independent identically distributed random points in $\T^d$ with common law the Lebesgue measure $\leb^d$. For every $i\in\I$, let $\Ncal_{i}=\intn{0}{-1+\sum_{j=0}^\infty m_{i,j}}$ if $\sum_{j=0}^\infty m_{i,j}<\infty$ and let $\Ncal_{i}=\N_{0}$ otherwise. The points $X_{i,n}$, for $n\in\Ncal_{i}$, are thus uniformly and independently chosen on the torus. Moreover, they are intended to govern the positions of the nonvanishing coefficients of the wavelet series $R$ in the direction $i$. The corresponding scales are
\[
j_{i,n}=\min\left\{ j\in\N_{0} \:\Biggl|\: n<\sum_{j'=0}^j m_{i,j'} \right\}.
\]
For $i\in\I$ and $n\in\Ncal_{i}$, let $\lambda_{i,n}$ denote the unique dyadic cube of generation $j_{i,n}$ that contains $X_{i,n}$. Finally, let
\[
\Mu=\left\{ (i,\lambda)\in\I\times\Lambda \:\bigl|\: \exists n\in \Ncal_{i}\quad\lambda=\lambda_{i,n} \right\} \qquad\text{and}\qquad R=\sum_{(i,\lambda)\in\Mu} 2^{-\underline{h} \gene{\lambda}} \Psi^i_{\lambda}.
\]
It is then easy to check that with probability one, for every $i\in\I$ and every $j\in\N_{0}$ large enough, there are exactly $m_{i,j}$ dyadic cubes $\lambda$ with generation $j$ such that $(i,\lambda)\in\Mu$, that is, such that the wavelet coefficient of $R$ indexed by $i$ and $\lambda$ does not vanish.

The absolute values of the coefficients of $R$ are at most $2^{-\underline{h}\gene{\lambda}}$. Proposition~4 in~\cite{Jaffard:2004fh} then ensures that $R$ belongs to the H\"older space $C^{\underline{h}}(\T^d)$. In addition, note that the pointwise regularity properties of $R$ are trivial when $\sum_{i,j}m_{i,j}<\infty$. In this case, there are finitely many pairs $(i,j)\in\I\times\N$ enjoying $m_{i,j}>0$, so that $\Mu$ is finite. The process $R$ is thereby $C^\infty$ everywhere on $\T^d$. In consequence, even if the results of this section are valid when the preceding sum is finite, they are nontrivial only when it is infinite.

\subsection{Statement of the results}

S.~Jaffard~\cite{Jaffard:2000mk} determined the law of the spectrum of singularities of $R$ in the case where $m_{i,j}=\lfloor 2^{\eta j}\rfloor$ for all $i\in\I$ and $j\in\N_{0}$ and some fixed $\eta\in (0,d)$. In order to refine this result and extend it to any family $(m_{i,j})_{(i,j)\in\I\times\N}$, we shall investigate the size and large intersection properties of
\[
E_{h}=\left\{ x\in\T^d \:\bigl|\: h_{R}(x)=h \right\} \qquad\text{and}\qquad \tilde E_{h}=\left\{ x\in\T^d \:\bigl|\: h_{R}(x)\leq h \right\}
\]
for every $h\in [0,\infty]$. Here, $h_{R}(x)$ denotes the H\"older exponent of $R$ at $x\in\T^d$. This exponent is defined as the supremum of all positive $h$ such that
\[
\dist(x',x)\leq\delta \qquad\Longrightarrow\qquad |R(x')-P(x'-x)|\leq c\ \dist(x',x)^h
\]
for all $x'\in\T^d$, some positive reals $c$ and $\delta$ and some function $P$ on $\T^d$ such that $P\circ\phi$ coincides with a polynomial in a neighborhood of zero, see~\cite{Jaffard:2004fh}.

For every $g\in\jauge_{d}$ and every nonempty open subset $V$ of $\T^d$, let $\grint^g(V)$ be the collection of all subsets $F$ of $\T^d$ such that $\phi^{-1}(F)$ belongs to the class $\grint^g(\phi^{-1}(V))$ defined in Section~\ref{grintlevyintro}. Theorem~\ref{GRINTSTABLE} then implies that the class $\grint^g(V)$ is closed under countable intersections and that each set in $\grint^g(V)$ has infinite Hausdorff $\overline{g}$-measure in every nonempty open subset of $V$ for every gauge function $\overline{g}\in\jauge_{d}$ enjoying $\overline{g}\prec g$. Furthermore, for every gauge function $g\in\jauge$, let
\[
h_{g}=\inf\left\{ h\in (0,\infty) \:\biggl|\: \sum_{i,j} m_{i,j}\ g_{d}(2^{-\underline{h}j/h})=\infty \right\}.
\]
Let $\overline{h}=h_{\Id^d}$. It is straightforward to check that $\overline{h}\geq\underline{h}$ and that $h_{\Id^s}=\overline{h}s/d$ for all $s\in (0,d]$. More generally, if $g_{d}=0$ then $h_{g}=\infty$, else $h_{g}\leq\overline{h}$.

The size and large intersection properties of the sets $E_{h}$ et $\tilde E_{h}$ are described by the following theorem, which can be thought of as the analog of Theorem~\ref{grintisoholdlevy}.

\begin{thm}\label{grintisoholdlac}
With probability one, $E_{h}=\emptyset$ for all $h\in [0,\underline{h})\cup(\overline{h},\infty]$ and $E_{\overline{h}}$ has full Lebesgue measure in $\T^d$. Moreover, let $u\in\jauge_{d}$ with $h_{u}<\infty$. With probability one, for every $v\in\jaugew_{1}$, every $h\in [\underline{h},\overline{h})$ and every nonempty open subset $V$ of $\T^d$,
\[
\hau^{u\circ v}(E_{h}\cap V)=\begin{cases}
0 & \text{si }h<h_{u\circ v} \\
\infty & \text{si }h_{u\circ v}\leq h\leq h_{u}
\end{cases}
\quad\text{and}\quad
\hau^{u\circ v}(\tilde E_{h}\cap V)=\begin{cases}
0 & \text{si }h<h_{u\circ v} \\
\infty & \text{si }h\geq h_{u\circ v}
\end{cases}
\]
and the set $\tilde E_{h}$ belongs to the class $\grint^{u\circ v}(V)$ if and only if $h\geq h_{u\circ v}$.
\end{thm}

\begin{rems}
This theorem yields an analog of Theorems~\ref{grintisoholdlevy1} and~\ref{grintisoholdlevy3}(\ref{grintisoholdlevy31}). To be specific, for any fixed gauge function $g\in\jauge$, with probability one, for every $h\in [\underline{h},\overline{h})$ and every nonempty open subset $V$ of $\T^d$,
\[
\hau^g(E_{h}\cap V)=\begin{cases}
0 & \text{if }h<h_{g} \\
\infty & \text{if }h=h_{g}
\end{cases}
\qquad\text{and}\qquad
\hau^g(\tilde E_{h}\cap V)=\begin{cases}
0 & \text{if }h<h_{g} \\
\infty & \text{if }h\geq h_{g}
\end{cases}
\]
and when $g\in\jauge_{d}$, the set $\tilde E_{h}$ belongs to $\grint^g(V)$ if and only if $h\geq h_{g}$. In particular, for $h\in [\underline{h},\overline{h})$, the set $E_{h}$ almost surely has infinite Hausdorff $g$-measure in every nonempty open subset of $\T^d$ if $g$ denotes a fixed gauge function in $\jauge$ such that $h_{g}=h$. However, such a function need not always exist. For instance, if $m_{i,j}=1$ for all $i\in\I$ and $j\in\N_{0}$, it is easy to verify that $h_{g}\in\{0,\infty\}$ for any $g\in\jauge$.

Theorem~\ref{grintisoholdlac} leads to an analog of Theorems~\ref{grintisoholdlevy2} and~\ref{grintisoholdlevy3}(\ref{grintisoholdlevy32}) as well. Indeed, if $\overline{h}<\infty$, then with probability one, for every $g\in\jauge$ with $g^{1/d}\in\jaugew$, every $h\in [\underline{h},\overline{h})$ and every nonempty open subset $V$ of $\T^d$,
\[
\hau^g(E_{h}\cap V)=\hau^g(\tilde E_{h}\cap V)=\begin{cases}
0 & \text{if }h<h_{g} \\
\infty & \text{if }h\geq h_{g}
\end{cases}
\]
and when $g\in\jauge_{d}$, the set $\tilde E_{h}$ belongs to $\grint^g(V)$ if and only if $h\geq h_{g}$. To prove this result, it suffices to apply Theorem~\ref{grintisoholdlac} with $u=\Id^d$ and $v={g_{d}}^{1/d}\in\jaugew_{1}$. As a consequence, with probability one, for all $h\in [\underline{h},\overline{h})$ and $s\in (0,d]$, the set $E_{h}$ has infinite (resp. zero) $s$-dimensional Hausdorff measure if $s\leq d h/\overline{h}$ (resp. $s>d h/\overline{h}$). The spectrum of singularities $h\mapsto\dim E_{h}$ of almost every sample path of the wavelet series $R$ is therefore given by
\[
\forall h\in [0,\infty] \qquad \dim E_{h}=\begin{cases}
dh/\overline{h} & \text{if }h\in [\underline{h},\overline{h}] \\
-\infty & \text{else}
\end{cases}.
\]
\end{rems}

We also study the size and large intersection properties of the set of all points at which a given function $w$ cannot be a modulus of continuity of the wavelet series $R$. Recall that the set $\modc$ is defined in Section~\ref{grintlevyintro}. As $R\in C^{\underline{h}}(\T^d)$, we may restrict ourselves to assuming that $w$ belongs to the set
\[
\modc_{\underline{h}}=\left\{ w\in\modc \:\bigl|\: w(\delta)={\rm o}(\delta^{\underline{h}}) \text{ as } \delta\to 0 \right\}.
\]
A function $w\in\modc_{\underline{h}}$ is a modulus of continuity of $R$ at a point $x\in\T^d$ if there are two positive reals $c$ and $\delta$, along with a function $P$ on $\T^d$ such that $P\circ\phi$ coincides with a polynomial in a neighborhood of zero and
\[
\dist(x',x)\leq\delta \qquad\Longrightarrow\qquad |R(x')-P(x'-x)|\leq c\ w(\dist(x',x))
\]
for all $x'\in\T^d$. Given $w\in\modc_{\underline{h}}$, let $F_{w}$ denote the set of all points at which $w$ is not a modulus of continuity of the wavelet series $R$. The following result can be seen as the analog of Theorem~\ref{grintmodclevy}.

\begin{thm}\label{grintmodclac}
Let $g\in\jauge$ and $w\in\modc_{\underline{h}}$ with $\sum_{i,j} m_{i,j} g_{d}(w^{-1}(2^{-\underline{h}j}))=\infty$. Then $F_{w}$ almost surely contains a set of the class $\grint^{g_{d}}(\T^d)$. Moreover, with probability one, $\hau^g(F_{w}\cap V)=\hau^g(V)$ for every open subset $V$ of $\T^d$.
\end{thm}

The remainder of this section is devoted to the proof of Theorems~\ref{grintisoholdlac} and~\ref{grintmodclac}. The methods which come into play are those introduced in Sections~\ref{grintlevyhold} and~\ref{grintlevymodc}.

\subsection{Proof of Theorem~\ref{grintisoholdlac}}

We begin by introducing some further notations and establishing two preliminary lemmas. The H\"older exponent of $R$ at a given point $x\in\T^d$ depends on how its nonvanishing wavelet coefficients are located around $x$. More precisely, for every continuous nondecreasing function $\ph:[0,\infty)\to\R$ enjoying $\ph(0)=0$, let
\[
L_{\ph}=\left\{ x\in\T^d \:\bigl|\: \dist(x,x_{\lambda})<\ph(2^{-\gene{\lambda}})\text{ for infinitely many }(i,\lambda)\in \Mu \right\}.
\]
Note that $\alpha\mapsto L_{\Id^{\underline{h}/\alpha}}$ is nondecreasing. The sets $L_{\Id^{\underline{h}/\alpha}}$ are related with the singularity sets $E_{h}$ and $\tilde E_{h}$, as shown by the following lemma, which can be regarded as the analog of Proposition~\ref{locholdlevy}.

\begin{lem}\label{locholdlac}
With probability one, for every $h\in [0,\underline{h})\cup(\overline{h},\infty]$, the set $E_{h}$ is empty and, for every $h\in [\underline{h},\overline{h}]$,
\[
\tilde E_{h}=\bigcap_{h<\alpha\leq\overline{h}} L_{\Id^{\underline{h}/\alpha}} \qquad\text{and}\qquad E_{h}=\tilde E_{h}\backslash\bigcup_{\underline{h}<\alpha<h} L_{\Id^{\underline{h}/\alpha}}.
\]
\end{lem}

\begin{proof}
The set $E_{h}$ is empty for all $h\in [0,\underline{h})$ because $R\in C^{\underline{h}}(\T^d)$. In addition, let $\alpha\in (\underline{h},\infty)$ and $x\in\T^d$. Proposition~1.3 in~\cite{Jaffard:1996kt} ensures that if $x\in L_{\Id^{\underline{h}/\alpha}}$ then $h_{R}(x)\leq\alpha$, else $h_{R}(x)\geq\alpha$. It follows that for any $h\in [\underline{h},\infty]$,
\[
\tilde E_{h}=\bigcap_{\alpha>h} L_{\Id^{\underline{h}/\alpha}} \qquad\text{and}\qquad E_{h}=\tilde E_{h}\backslash\bigcup_{\underline{h}<\alpha<h} L_{\Id^{\underline{h}/\alpha}}.
\]
To conclude, observe that for all $\alpha\in (\overline{h},\infty)$, the set $L_{\Id^{\underline{h}/\alpha}}$ contains the set of all $x\in\T^d$ such that $\dist(x,X_{i,n})<2^{-1-\underline{h}j_{i,n}/\alpha}$ for infinitely many $i\in\I$ and $n\in\Ncal_{i}$, which is almost surely equal to $\T^d$, by Proposition~9 in~\cite[Chapter~11]{Kahane:1985gc}.
\end{proof}

The following result is analogous to Lemma~\ref{sizelph}.

\begin{lem}\label{sizelphlac}
Consider a continuous nondecreasing function $\ph:[0,\infty)\to\R$ with $\ph(0)=0$. If the sum $\sum_{i,j} m_{i,j}\ph(2^{-j})$ converges, then with probability one,
\[
\sum_{(i,\lambda)\in\Mu} \ph(2^{-\gene{\lambda}})<\infty.
\]
If it diverges and $r^d={\rm o}(\ph(r))$ as $r\to 0$, then the set $L_{\ph^{1/d}}$ almost surely has full Lebesgue measure in $\T^d$.
\end{lem}

\begin{proof}
Assume that $\sum_{i,j} m_{i,j}\ph(2^{-j})<\infty$ and let $(i,\lambda)\in\I\times\Lambda$. Then, $(i,\lambda)\in\Mu$ if and only if $X_{i,n}\in\lambda$ for some $n\in\intn{m_{i,0}+\ldots+m_{i,\gene{\lambda}-1}}{m_{i,0}+\ldots+m_{i,\gene{\lambda}}-1}$. Thus, $(i,\lambda)\in\Mu$ with probability $1-(1-2^{-d\gene{\lambda}})^{m_{i,\gene{\lambda}}}\leq 2^{-d\gene{\lambda}}m_{i,\gene{\lambda}}$. Therefore,
\[
\esp\left[\sum_{(i,\lambda)\in\Mu}\ph(2^{-\gene{\lambda}})\right]\leq\sum_{(i,\lambda)\in\I\times\Lambda} 2^{-d\gene{\lambda}}m_{i,\gene{\lambda}}\ph(2^{-\gene{\lambda}})=\sum_{(i,j)\in\I\times\N} m_{i,j}\ph(2^{-j})<\infty,
\]
so that the sum $\sum_{(i,\lambda)\in\Mu}\ph(2^{-\gene{\lambda}})$ is finite with probability one.

Conversely, suppose that $\sum_{i,j} m_{i,j}\ph(2^{-j})=\infty$ and $r^d={\rm o}(\ph(r))$. Hence, $L_{\ph^{1/d}}$ contains the set $\tilde L_{\ph^{1/d}}$ of all $x\in\T^d$ such that $\dist(x,X_{i,n})<\ph(2^{-j_{i,n}})^{1/d}/2$ for infinitely many $i\in\I$ and $n\in\Ncal_{i}$. Let $x\in\T^d$ and assume that $x\not\in\tilde L_{\ph^{1/d}}$. There is a positive integer $j_{0}$ such that $\dist(X_{i,n},x)\geq\ph(2^{-j})^{1/d}/2$ for every $i\in\I$, every integer $j\geq j_{0}$ and every integer $n$ between $m_{i,0}+\ldots+m_{i,j-1}$ and $m_{i,0}+\ldots+m_{i,j}-1$. This occurs with probability at most
\[
\prod_{i\in\I \atop j\geq j_{0}} \left(1-\kappa\,\ph(2^{-j})\right)^{m_{i,j}}\leq\exp\left(-\kappa \sum_{i\in\I \atop j\geq j_{0}} m_{i,j}\ph(2^{-j})\right)
\]
where $\kappa>0$ only depends on the norm $\R^d$ is endowed with. Note that the right-hand side vanishes. As a consequence, every point $x\in\T^d$ almost surely belongs to $L_{\ph^{1/d}}$. We conclude thanks to Fubini's theorem.
\end{proof}

Let us establish Theorem~\ref{grintisoholdlac}. For the sake of clarity, we divide its statement into six propositions, namely, Propositions~\ref{prpEhemptyfullmeas} to~\ref{hauuvehinftylac}, which we now announce and prove.

\begin{prp}\label{prpEhemptyfullmeas}
With probability one, the set $E_{h}$ is empty for all $h\in [0,\underline{h})\cup(\overline{h},\infty]$ and the set $E_{\overline{h}}$ has full Lebesgue measure in $\T^d$.
\end{prp}

\begin{proof}
To begin with, Lemma~\ref{locholdlac} implies that with probability one, the set $E_{h}$ is empty for all $h\in [0,\underline{h})\cup(\overline{h},\infty]$. This result also ensures that with probability one,
\[
E_{\overline{h}}=\T^d\backslash\bigcup_{\underline{h}<\alpha<\overline{h}} L_{\Id^{\underline{h}/\alpha}}.
\]
Let $\alpha\in(\underline{h},\overline{h})$. For each $j\in\N$, the set $L_{\Id^{\underline{h}/\alpha}}$ is covered by the open balls with center $x_{\lambda}$ and radius $2^{-\underline{h}\gene{\lambda}/\alpha}$ for all $(i,\lambda)\in\Mu$ enjoying $\gene{\lambda}\geq j$. This covering yields
\[
\leb^d(L_{\Id^{\underline{h}/\alpha}})\leq\kappa\sum_{(i,\lambda)\in\Mu \atop \gene{\lambda}\geq j} 2^{-d\underline{h}\gene{\lambda}/\alpha}
\]
where $\kappa>0$ only depends on the norm $\R^d$ is endowed with. Lemma~\ref{sizelphlac} implies that the right-hand side tends to zero as $j\to\infty$. Since $\alpha\mapsto L_{\Id^{\underline{h}/\alpha}}$ is nondecreasing, it follows that with probability one, for every $\alpha\in (\underline{h},\overline{h})$, the set $L_{\Id^{\underline{h}/\alpha}}$ has Lebesgue measure zero. Thus, $E_{\overline{h}}$ almost surely has full Lebesgue measure in $\T^d$.
\end{proof}

Let $u\in\jauge_{d}$ with $h_{u}<\infty$ and let $\tilde u$ denote a continuous nondecreasing function defined on $[0,\infty)$ that coincides with $u$ in a neighborhood of the origin. Since $h_{u}=h_{\tilde u}$, Lemma~\ref{sizelphlac} implies that with probability one, for every real $h>\underline{h}$,
\begin{equation}\label{sizelphlaceq}\begin{split}
& h<h_{u} \qquad\Longrightarrow\qquad \sum\limits_{(i,\lambda)\in\Mu} \tilde u(2^{-\underline{h}\gene{\lambda}/h})<\infty \\
\text{and}\qquad & h>h_{u} \qquad\Longrightarrow\qquad \leb^d(L_{(\tilde u\circ\Id^{\underline{h}/h})^{1/d}})=1.
\end{split}\end{equation}
Thus~(\ref{sizelphlaceq}) and the statement of Lemma~\ref{locholdlac} hold with probability one. Assume that the corresponding event occurs, let $v\in\jaugew_{1}$ and let $\tilde v$ denote a continuous nondecreasing function defined on $[0,\infty)$ that coincides with $v$ in a neighborhood of the origin. The limit
\[
\gamma_{v}=\lim_{r\to 0}\frac{\log v(r)}{\log r}
\]
exists and belongs to $[0,1]$. Moreover, the functions $u\circ v$ and $\tilde u\circ\tilde v$ coincide in a neighborhood of zero and $h_{\tilde u\circ\tilde v}=h_{u\circ v}=h_{u}\gamma_{v}$.

Let $h\in [\underline{h},\overline{h})$ and let $V$ be a nonempty open subset of $\T^d$.

\begin{prp}\label{hauuveh0lachau}
If  $h<h_{u\circ v}$, then $\hau^{u\circ v}(E_{h}\cap V)=\hau^{u\circ v}(\tilde E_{h}\cap V)=0$.
\end{prp}

\begin{proof}
Because of Lemma~\ref{locholdlac}, we have $E_{h}\subseteq\tilde E_{h}\subseteq L_{\Id^{\underline{h}/\alpha}}$ for any $\alpha\in (h,h_{u\circ v})$. Hence it suffices to show that $L_{\Id^{\underline{h}/\alpha}}$ has zero Hausdorff $u\circ v$-measure. To this end, note that $\underline{h}<\alpha/(\gamma_{v}-\eta)<h_{u}$ for some $\eta\in (0,\gamma_{v})$. Besides, for every $\eps>0$ small enough, $u\circ v(r)\leq\tilde u(r^{\gamma_{v}-\eta})$ for all $r\in [0,\eps]$ and the set $L_{\Id^{\underline{h}/\alpha}}$ is covered by the open balls with center $x_{\lambda}$ and radius $2^{-\underline{h}\gene{\lambda}/\alpha}$ for all $(i,\lambda)\in\Mu$ satisfying $2^{1-\underline{h}\gene{\lambda}/\alpha}<\eps$. This covering yields
\[
\hau^{u\circ v}_{\eps}(L_{\Id^{\underline{h}/\alpha}})\leq 2^d\sum_{(i,\lambda)\in\Mu \atop 2^{1-\underline{h}\gene{\lambda}/\alpha}<\eps} \tilde u(2^{-\underline{h}\gene{\lambda}(\gamma_{v}-\eta)/\alpha}).
\]
As shown by~(\ref{sizelphlaceq}), the right-hand side tends to zero as $\eps\to 0$. Hence $L_{\Id^{\underline{h}/\alpha}}$ has zero $u\circ v$-measure.
\end{proof}

\begin{prp}\label{hauuveh0lacgrint}
If $h<h_{u\circ v}$, then $\tilde E_{h}\not\in\grint^{u\circ v}(V)$.
\end{prp}

\begin{proof}
There exists $\overline{u}\in\jauge_{d}$ with $\overline{u}\prec u$ and $h_{\overline{u}}=h_{u}<\infty$. Proposition~\ref{hauuveh0lachau} still holds with $\overline{u}$ instead of $u$, so that $\hau^{\overline{u}\circ v}(\tilde E_{h}\cap V)=0$. Proposition~\ref{hauuveh0lacgrint} is then a direct consequence of Theorem~\ref{GRINTSTABLE}, together with the observation that $\overline{u}\circ v\prec u\circ v$.
\end{proof}

\begin{prp}\label{hauuvtildeehinftylacgrint}
If $h\geq h_{u\circ v}$, then $\tilde E_{h}\in\grint^{u\circ v}(V)$.
\end{prp}

\begin{proof}
Let $\alpha\in(h,\overline{h}]$. Note that $\alpha/(\gamma_{v}+\eps)>\max(h_{u},\underline{h})$ for $\eps>0$ small enough. Moreover, the set $L_{(\tilde u\circ\tilde v\circ\Id^{\underline{h}/\alpha})^{1/d}}$, containing $L_{(\tilde u\circ\Id^{\underline{h}(\gamma_{v}+\eps)/\alpha})^{1/d}}$, has full Lebesgue measure in $\T^d$, owing to~(\ref{sizelphlaceq}). Hence, the family $(p+\dot x_{\lambda},\tilde u\circ\tilde v(2^{-\underline{h}\gene{\lambda}/\alpha})^{1/d})_{(i,\lambda,p)\in\Mu\times\Z^d}$ is a homogeneous ubiquitous system in $\R^d$. Here, $\dot x_{\lambda}$ denotes the unique element of $\phi^{-1}(\{x_{\lambda}\})\cap [0,1)^d$. Theorem~\ref{UBIQUITYHOM} then implies that $L_{\Id^{\underline{h}/\alpha}}\in\grint^{u\circ v}(\T^d)$. Because of Lemma~\ref{locholdlac}, the set $\tilde E_{h}$ can be written as a countable intersection of the sets $L_{\Id^{\underline{h}/\alpha}}$ for $\alpha\in (h,\overline{h}]$. As the class $\grint^{u\circ v}(\T^d)$ is closed under countable intersections, it necessarily contains $\tilde E_{h}$. Proposition~\ref{hauuvtildeehinftylacgrint} follows.
\end{proof}

\begin{prp}\label{hauuvtildeehinftylachau}
If $h\geq h_{u\circ v}$, then $\hau^{u\circ v}(\tilde E_{h}\cap V)=\infty$.
\end{prp}

\begin{proof}
Let us first assume that $u\not\prec\Id^d$. Thus $\overline{h}\leq h_{u}<\infty$ and $\gamma_{v}<1$, so that $v\prec\Id$. Proposition~\ref{hauuvtildeehinftylacgrint} with $\Id^d$ instead of $u$ and $\underline{v}:r\mapsto v(r)/\log(v(r)/r)$ instead of $v$ shows that $\tilde E_{h}\in\grint^{\underline{v}^d}(V)$. As $v^d\prec\underline{v}^d$, this set has infinite Hausdorff $v^d$-measure in $V$ and the result follows. Conversely, if $u\prec\Id^d$, there exists a gauge function $\underline{u}\in\jauge_{d}$ with $u\prec\underline{u}$ and $h_{\underline{u}}\leq h_{u}$. Proposition~\ref{hauuvtildeehinftylacgrint} with $\underline{u}$ rather than $u$ implies that $\tilde E_{h}\in\grint^{\underline{u}\circ v}(V)$. The result thereby follows from the fact that $u\circ v\prec\underline{u}\circ v$.
\end{proof}

\begin{prp}\label{hauuvehinftylac}
If $h\in [h_{u\circ v},h_{u}]$, then $\hau^{u\circ v}(E_{h}\cap V)=\infty$.
\end{prp}

\begin{proof}
Let us suppose that $h=h_{u\circ v}$ and let $\alpha\in (\underline{h},h)$. There exists $\eta\in (0,\gamma_{v})$ such that $\underline{h}<\alpha/(\gamma_{v}-\eta)<h_{u}$. For $\eps>0$ small enough and every $r\in [0,\eps]$, we have $u\circ v(r)\leq\tilde u(r^{\gamma_{v}-\eta})$ and the set $L_{\Id^{\underline{h}/\alpha}}$ is covered by the open balls with center $x_{\lambda}$ and radius $2^{-\underline{h}\gene{\lambda}/\alpha}$ for all $(i,\lambda)\in\Mu$ enjoying $2^{1-\underline{h}\gene{\lambda}/\alpha}<\eps$. This covering yields
\[
\hau^{u\circ v}_{\eps}(L_{\Id^{\underline{h}/\alpha}})\leq 2^d\sum_{(i,\lambda)\in\Mu \atop 2^{1-\underline{h}\gene{\lambda}/\alpha}<\eps} \tilde u(2^{-\underline{h}\gene{\lambda}(\gamma_{v}-\eta)/\alpha}).
\]
The right-hand side tends to zero as $\eps\to 0$ by virtue of~(\ref{sizelphlaceq}), so that $L_{\Id^{\underline{h}/\alpha}}$ has zero $u\circ v$-measure. Lemma~\ref{locholdlac} and Proposition~\ref{hauuvtildeehinftylachau} then lead to Proposition~\ref{hauuvehinftylac} for $h=h_{u\circ v}$. Conversely, if $h>h_{u\circ v}$, we may rewrite what precedes with $w:r\mapsto r^{h/h_{u}}$ instead of $v$ in order to obtain $\hau^{u\circ w}(E_{h}\cap V)=\infty$. The result finally follows from the observation that $u\circ v(r)\geq u\circ w(r)$ for $r\geq 0$ small enough.
\end{proof}

\subsection{Proof of Theorem~\ref{grintmodclac}}

Consider two functions $g\in\jauge$ and $w\in\modc_{\underline{h}}$ with $\sum_{i,j}m_{i,j} g_{d}(w^{-1}(2^{-\underline{h}j}))=\infty$. Let $\tilde g$ be a continuous nondecreasing function defined on $[0,\infty)$ that coincides with $g_{d}$ in a neighborhood of zero and let $\tilde w$ denote a continuous increasing function defined on $[0,\infty)$ that tends to infinity at infinity and coincides with $w$ in a neighborhood of the origin. Note that $\kappa\,\tilde w(\delta)\leq\tilde w(2\delta)$ for every $\delta\in [0,\delta_{0}]$ and some $\kappa>1$ and $\delta_{0}>0$. For each $q\in\N$, let $\ph_{q}:r\mapsto\tilde w^{-1}(r^{\underline{h}}/\kappa^q)$ and let
\[
\tilde F_{w}=\bigcap_{q=0}^\infty\downarrow L_{\ph_{q}}.
\]
Theorem~\ref{grintmodclac} is a direct consequence of the two following lemmas.

\begin{lem}
With probability one, $\tilde F_{w}\in\grint^{g_{d}}(\T^d)$ and $\hau^g(\tilde F_{w}\cap V)=\hau^g(V)$ for every open subset $V$ of $\T^d$.
\end{lem}

\begin{proof}
First note that $\sum_{i,j} m_{i,j}\tilde g(\ph_{q}(2^{-j}))=\infty$, because $\tilde g(\ph_{q}(r))\geq\tilde g(\tilde w^{-1}(r^{\underline{h}}))$ for $r\geq 0$ small enough. Suppose that $g_{d}\not\prec\Id^d$. The measure $\hau^g$ coincides up to a multiplicative constant with the Lebesgue measure on the Borel subsets of $\T^d$ and $\sum_{i,j} m_{i,j}\ph_{q}(2^{-j})^d=\infty$, so that $L_{\ph_{q}}$ almost surely has full Lebesgue measure in $\T^d$ by Lemma~\ref{sizelphlac}. It follows that with probability one, $\hau^g(\tilde F_{w}\cap V)=\hau^g(V)$ for every open $V\subseteq\T^d$. In the general case, Lemma~\ref{sizelphlac} and Theorem~\ref{UBIQUITYHOM} imply that $L_{\ph_{q}}$ almost surely belongs to $\grint^{g_{d}}(\T^d)$. As this class is closed under countable intersections, it almost surely contains $\tilde F_{w}$. In addition, if $g_{d}\prec\Id^d$, there is a gauge function $\underline{g}\in\jauge_{d}$ with $g_{d}\prec\underline{g}$ and $\sum_{i,j}m_{i,j} \underline{g}(w^{-1}(2^{-\underline{h}j}))=\infty$. We may use $\underline{g}$ instead of $g_{d}$ above in order to prove that $\tilde F_{w}\in\grint^{\underline{g}}(\T^d)$ with probability one. This yields $\hau^{g_{d}}(\tilde F_{w}\cap V)=\infty=\hau^{g_{d}}(V)$ for every nonempty open subset $V$ of $\T^d$. Hence, with probability one, $\hau^g(\tilde F_{w}\cap V)=\hau^g(V)$ for every open $V\subseteq\T^d$.
\end{proof}

\begin{lem}
We have $\tilde F_{w}\subseteq F_{w}$.
\end{lem}

\begin{proof}
Let $x\in\tilde F_{w}$ and suppose that $x\not\in F_{w}$. Thus, $w$ is a modulus of continuity of $R$ at $x$ and Proposition~3 in~\cite{Jaffard:2000mk} implies that there exists a real $c>0$ such that
\[
2^{-\underline{h}\gene{\lambda}}\un_{\{(i,\lambda)\in\Mu\}}\leq c\,(w(2^{-\gene{\lambda}})+w(\dist(x,x_{\lambda})))
\]
for every $i\in\I$ and every $\lambda\in\Lambda$ with $\gene{\lambda}$ large enough and $\dist(x,x_{\lambda})$ small enough. Let $q$ be large enough to ensure that $\kappa^q>2c$. As $x\in F_{\ph_{q}}$, there are infinitely many $(i,\lambda)\in\Mu$ with $\dist(x,x_{\lambda})<\tilde w^{-1}(2^{\underline{h}\gene{\lambda}}/\kappa^q)$. Hence, for $\gene{\lambda}$ large enough,
\[
2^{-\underline{h}\gene{\lambda}}\leq c\,(w(2^{-\gene{\lambda}})+w(\dist(x,x_{\lambda})))\leq \frac{2c}{\kappa^q}2^{-\underline{h}\gene{\lambda}}
\]
which is a contradiction.
\end{proof}

\end{document}